\documentclass[a4paper,12pt,dvipdfmx]{article}

\usepackage{amsthm}
\usepackage{amsmath,amssymb,latexsym,amsfonts,mathrsfs}
\usepackage{graphicx}

\usepackage{color}

\usepackage{amsrefs}
\usepackage{mathtools}
\usepackage{fancyhdr}
\usepackage[top=30truemm, bottom=30truemm, left=25truemm, right=25truemm]{geometry}

\allowdisplaybreaks[4]

\newcommand{\ep}{\varepsilon}

\newcommand{\SCR}[1]{{\mathscr #1}}
\newcommand{\CAL}[1]{{\cal #1}}
\newcommand{\J}[1]{\left\langle #1 \right\rangle}
\newcommand{\D}[1]{{\mathscr D}( #1 )}

\theoremstyle{definition}
\newtheorem{Thm}{{\bf Theorem}}[section]
\newtheorem{Lem}[]{{\bf Lemma}}[section]
\newtheorem{Prop}[]{{\bf Proposition}}[section]
\newtheorem{Cor}[]{{\bf Corollary}}[section]
\newtheorem{Assu}[]{{\bf Assumption}}[section]
\newtheorem{Def}[]{{\bf Definition}}[section]
\newtheorem{Rem}[]{{\bf Remark}}[section]

\usepackage{diffcoeff}
\usepackage[italicdiff]{physics}

\numberwithin{equation}{section}

\usepackage{bm}
\usepackage{hyperref} 
\hypersetup{
    colorlinks=true, 
    linktoc= page,   
citecolor=blue
}

\def\({\left(}
\def\){\right)}
\def\<{\left\langle}
\def\>{\right\rangle}

\newcommand{\C}{{\bf C}}
\newcommand{\R}{{\bf R}}
\newcommand{\Z}{{\bf Z}}
\newcommand{\N}{{\bf N}}

\renewcommand{\Re}{\operatorname{Re}}

\newcommand{\relmiddle}[1]{\mathrel{}\middle#1\mathrel{}}

\begin{document}
\begin{flushleft}
{\Large \bf Modified scattering operator for nonlinear Schr\"odinger equations with time-decaying harmonic potentials}
\end{flushleft}

\begin{flushleft}
{\large Masaki KAWAMOTO}\\
{Department of Engineering for Production, Graduate School of Science and Engineering, Ehime University, 3 Bunkyo-cho Matsuyama, Ehime, 790-0826, Japan }\\
Email: {kawamoto.masaki.zs@ehime-u.ac.jp}
\end{flushleft}
\begin{flushleft}
{\large Hayato MIYAZAKI}\\
{Teacher Training Courses, Faculty of Education, Kagawa University, Takamatsu, Kagawa 760-8522, Japan} \\
{E-mail: miyazaki.hayato@kagawa-u.ac.jp}
\end{flushleft}

\begin{abstract}
This paper is concerned with nonlinear Schr\"odinger equations with a time-decaying harmonic potential.
The nonlinearity is gauge-invariant of the long-range critical order. 
In \cite{KM21} and \cite{Ka21}, it is proved that the equation admits a nontrivial solution that behaves like a free solution with a logarithmic phase correction in the frameworks of both the final state problem and the initial value problem.
Furthermore, a modified scattering operator has been established in the case without the potential in \cite{HN06}. 
In this paper, we construct a modified scattering operator for our equation by utilizing a generator of the Galilean transformation.
Moreover, we remove a restriction for the coefficient of the potential which is required in \cite{Ka21}.
\end{abstract}

\begin{flushleft}
2020 \textit{Mathematics Subject Classification}: Primary: 35Q55, Secondly: 35B40, 35P25.

\textit{Key words and phrases}: nonlinear Schr\"odinger equation; time-decaying harmonic potential; long-range scattering; modified scattering operator
\end{flushleft}



\section{Introduction}
\subsection{Setting and Backgounrd}
In this paper, we consider the nonlinear Schr\"{o}dinger equation 
\begin{align}\label{eq1} \tag{NLS}
	i \partial_t u - H(t) u =  F(u),
\end{align}
where $(t,x) \in \R \times \R^d$, $d = 1,2,3$, and $u = u(t,x)$ is a complex-valued unknown function.
The operator $H(t)$ is defined by
\begin{align*}
	H (t) = - \frac{1}{2}\Delta + \sigma (t) \frac{|x|^2}{2}, \quad \sigma (t) \in L^{\infty}(\R; \R).
\end{align*}
The coefficient of the potential $\sigma (t)$ satisfies the time-decay condition $\lim_{t \to \infty} t^2\sigma(t) \eqqcolon \sigma_1 \in (- \infty, 1/4)$, more precisely see Assumption \ref{A1}. The nonlinearity $F$ is gauge invariant of degree $1+ 2/d(1-\lambda)$; that is, it is the form of 
\begin{align}
	F(u) = \eta |u|^{p_{c}} u, \quad p_c \coloneqq \frac{2}{d(1-\lambda)}, \quad 
	\eta \in \R,  \label{noho}
\end{align}
where $\lambda = (1- \sqrt{1-4 \sigma_1} ) /2 < 1/2$ with $-1/3 < \lambda$ if $d=3$. 

The aim of this paper is to deal with a scattering operator to \eqref{eq1}. As for \eqref{eq1} with $ \sigma(t) \equiv 0$, that is,
\begin{align}\label{eq:f}
	i \partial_t u + \frac12 \Delta u =  G(u), \quad G(u) = \eta |u|^p u ,
\end{align}
then \eqref{eq:f} admits a nontrivial solution that behaves like a free solution $U_{0}(t) u_{\pm}$
as $t \rightarrow \pm \infty$, when $p> 2/d$, where $U_{0}(t) = e^{it \Delta/2}$ and $u_{ \pm}$ is called the \emph{final date}. In this sense, we say the case $p > 2/d$ is the {\em short-range}. 
On the other hand, if $p \leq 2/d$, then \eqref{eq:f} has no solutions that behave like the free solution in $L^2$ (e.g., \cite{S74, B84, TY84}), hence we say the case $p \leq d/2$ is the {\em long-range}. Hence we can regard  the exponent $p_c' \coloneqq 2/d$ as the threshold in the sense of scattering problems.  
In particular, when $p= p_c'$, there is the solution that behaves like a free solution with a logarithmic phase correction $U_N(t) u_{\pm} := U_{0}(t) \mathcal{F}^{-1} \exp\(-i \eta \left|\mathcal{F} u_{\pm}\right|^{p_{c}} \log |t|\) \mathcal{F} u_{\pm}$
for a suitable given $u_{\pm}$ as $t \to \pm \infty$ (cf. \cite{O91, GO93, HN98}), where $ \mathcal{F}$ is the usual Fourier transform.
A scattering problem for nonlinear Schr\"odinger equations has been considered in two directions.
One direction is the final state problem. 
We attempt to construct a solution $u$ to \eqref{eq:f} that behaves like a prescribed asymptotics $U_{N}(t) u_{-}$ as $t \rightarrow  - \infty$ in the problem.
The second is the initial value problem that is dealt with a global-in-time solution and asymptotics $U_{N}(t) u_{+}$ for large time, specified by given initial data. If we find the solution under the final state problem, then the wave operator $ \mathcal{W}_{-}$ which is a map from $u_{-}$ into $u(0)$ can be defined. 
On the other hand, we can define the inverse wave operator $\mathcal{W}_{+}^{-1}$ from $u(0)$ into $u_{+}$, once the initial value problem is solved.
Furthermore, if the range of $ \mathcal{W}_{-}$ included in the domain of $ \mathcal{W}_{+}^{-1}$, then the scattering operator $ \mathcal{S} \coloneqq \mathcal{W}_{+}^{-1} \mathcal{W}_{-}$ can be defined.

The scattering operator on the critical case $p = p_{c}'$ is called the \emph{modified scattering operator}.
In \cite{Ca01}, when $d=1$, Carles constructs a modified scattering operator for \eqref{eq:f} on $H^{0, 3} \cap H^{1, 2}$ if $\norm{u_{-}}_{H^{1} \cap H^{0,1}}$ is small.
Later on, the operator on the neighborhood centered at the origin in $H^{\alpha, 0} \cap H^{0, \alpha}$ has been defined in Hayashi and Naumkin \cite{HN06} under $d/2 < \alpha < \min(1+p_{c}', d)$. 
This situation such that the operator can be defined on the neighborhood centered at the origin is called \emph{low-energy scattering} (cf. \cite{Stra81a, Stra81b}).
This kind of modified scattering operator has been treated on other types of nonlinearity in some literature (e.g., \cite{Na02a, Na02b}).
In terms of \eqref{eq1}, Carles \cite{Ca03} proves the existence of the scattering operator, provided $\sigma(t) \equiv - \omega^{2}$ ($\omega>0$) and $F(u) = \eta |u|^{2 \sigma} u$, $\sigma >0$ when $d=1,2$ and $0 < \sigma < 2/(d-2)$ when $d \geq 3$.

Now we back to the case $\sigma(t) \neq 0$. Under the existence of $\sigma(t) |x| ^2$, the dispersive effect for free solution (\eqref{eq1} with $F\equiv 0$) becomes weaker than the case of $\sigma (t) \equiv 0$ and a threshold in the sense of scattering is changed from $p_c'$ to $p_c = 2/(d(1- \lambda))$, (see below). Hence the construction of the scattering operator for \eqref{eq1} is equivalent to that of the modified scattering operator. 
Kawamoto and Muramatsu \cite{KM21} study the initial value problem for \eqref{eq1} with $F(u) \sim \eta |u|^{p_c} u$, and they succeed to construct the inverse modified wave operator $\CAL{W}^{-1}_+$ from $u_0 \in H^{0,\alpha}$ to $u_+ \in L^2$. 
Also, Kawamoto \cite{Ka21} proves that the existence of modified wave operators $\CAL{W}^{}_{-}$ from $H^{0, \alpha}$ to $u_0 \in L^2$. Clearly, these results do NOT guarantee $\mathrm{Ran} (\CAL{W}_-) \subset \SCR{D}(\CAL{W}_+^{-1})$ and hence to construct the modified wave operator is difficult only with employing previous results. Hence the first aim of this paper is to bridge this gap. 
Moreover, a strong restriction for $\lambda$ has been imposed in \cite{KM21} because of a use of the time-weighted Strichartz estimate developed by \cite{KY18} and \cite{Ka20}. Hence the second aim of this paper is to remove the restriction. The modified scattering operator for \eqref{eq1} has not been discussed and so we attempt to define the operator in this paper under the some smallness conditions which are associated to the low energy condition in the case of $\sigma (t) =0$.

To state the reason why $p \leq p_c$ becomes long-range, we first let $U(t,s)$ be a propagator for $H (t)$, defined by the family of unitary operators $\{ U(t,s)\}_{(t,s) \in {\bf R}^2} $ on $L^2({\bf R}^n)$ such that for all $t,$ $s$,  $\tau \in {\bf R}$,
\begin{align*}
& i \partial _t U(t,s) = H(t) U(t,s), \quad i \partial _s U(t,s) = -U(t,s) H(s), \\
& U(t, \tau) U( \tau, s) = U(t,s), \quad U(s,s) = \mathrm{Id}_{L^2({\bf R}^n)}, \quad U(t,s) \D{H(s)} \subset \D{H(s) }
\end{align*}
hold on $\D{H(s)}$, i.e., $u_0(t) :=U(t,0) u_0$ solves $i \partial_t u_0(t) = H(t) u_0(t)$ with $u(0) = u_0 \in \D{H(0)}$. We here mention that  $U(t,s)$ has the following factorization formula 
\begin{align} \label{mdfm1}
	U(t,0) = \CAL{M} \left(  \frac{\zeta _2(t)}{\zeta _2 '(t)} \right) \CAL{D} (\zeta _2 (t)) \mathcal{F} \CAL{M} \left( \frac{\zeta _2 (t)}{ \zeta _1 (t)} \right),
\end{align}
where $\CAL{M}(t) = e^{\frac{i|x|^2}{2t}}$, $\left(\CAL{D}(t) \phi \right) (x) = (it)^{-\frac{d}{2}} \phi (x/t)$, $ \mathcal{F}$ is a usual Fourier transform and $\zeta_{j} (t)$, $j=1,2$ are solutions to 
\begin{align}
\zeta _j ''(t) + \sigma (t) \zeta _j (t) =0, \quad
\begin{cases}
\zeta _1 (0) = 1, \\
\zeta _1 ' (0) =0,
\end{cases}
\quad
\begin{cases}
\zeta _2 (0) = 0, \\
\zeta _2 ' (0) =1.
\end{cases}
	\label{ode:1}
\end{align}
Korotyaev \cite{Ko89} first discovered  similar factorization to \eqref{mdfm1} and Kawamoto and Muramatsu \cite{KM21} rewrites Korotyaev's factorization to \eqref{mdfm1} to deal with a nonlinear scattering problem on \eqref{eq1} motivated by Hayashi and Namunkin \cite{HN98} (see also \cite{Ca03, Ca11, Ka21}). This factorization formula leads to a kind of dispersive estimate
\begin{align}
	\norm{U(t,0) \phi}_{ \infty} \lesssim | \zeta_{2}(t)|^{-\frac{d}2} \norm{ \phi}_{1} .
	\label{disp}
\end{align}
The typical example of $\sigma(t)$ in our mind is 
\begin{align}
	\sigma (t) = 
	\begin{cases}
	\sigma_0 \quad (|t| < r_1), \\
	\sigma_1 |t|^{-2} \quad (|t| \geq r_1),
	\end{cases}
	\quad r_{1}, \sigma_{0} >0.
	\label{exsi:1}
\end{align} 
For this case, $\zeta_j(t)$, $|t| \geq r_1$ can be written as a linear combination of $|t|^{\lambda}$ and $|t|^{1- \lambda}$ with $\lambda = (1-\sqrt{1-4 \sigma_1})/2$. Namely, those solutions are
\begin{align*}
	\zeta_{1} (t) = b_{1} |t|^{1-\lambda} + b_{1}' |t|^{\lambda}, \quad 
	\zeta_{2} (t) = b_{2} |t|^{1-\lambda} + b_{2}' |t|^{\lambda} .
\end{align*}
Hence the time decay estimate becomes $\| U(t,0) \phi \|_{\infty} \lesssim |t|^{- \frac{d}2(1- \lambda)}$ holds whenever $b_2' \neq 0$,  and this decay is slower than the free case $\sigma(t) \equiv 0$. The fact changes a threshold that the asymptotics involves nonlinear effects coming from $G(u) = \eta |u|^{p} u$ with $p>0$.

\subsection{Main results}
In what follows, we assume the following implicit condition on $\sigma (t)$:
\begin{Assu} \label{A1}
There exist solutions $\zeta (t)$ to \eqref{ode:1} satisfying the following condition: 
There exist $r_0 \geq 1$ such that
\begin{enumerate}
\item[(A1)] there exists $\delta_{0} >0$ such that $|\zeta _2 (t)| \geq \delta_{0}$ for any $|t| \geq r_{0}$,
\item[(A2)] there exist $a_{j} \in {\bf R}$, $j=1,2$ with $a_{2} \neq 0$ and $\lambda < 1/2$ with $-1/3 < \lambda$ if $d=3$ such that
\begin{align}
	\abs{ \frac{\zeta_2 (t)}{|t|^{1-\lambda}} - a_{2}} +
	\abs{\frac{\zeta _1(t)}{\zeta _2 (t)} - \frac{a_{1}}{a_{2}}}
	\lesssim{}& 
	\begin{cases}
	|t|^{ -(1- 2 \lambda)} \quad ( \lambda \geq 0), \\
	|t|^{ -1} \quad ( \lambda < 0),
	\end{cases}
	\label{K13} \\
	\abs{ \frac{\zeta_2 (t)}{|t|^{1-\lambda}} - a_{2}} \leq \frac{|a_{2}|}{2}
	\label{K13-2}
\end{align}
for any $|t| \geq r_{0}$.
\end{enumerate}
\end{Assu}

We remark that the condition \eqref{K13} are required in order that $\mathcal{M}\left( \zeta _2 (t) / \zeta _1 (t) \right)$ is approximated $ \mathcal{M}(a_{2}/a_{1}) \eqqcolon \mathcal{M}_{+}$ as $t \rightarrow \pm \infty$. More precisely, we obtain
\begin{align}
	\abs{\mathcal{M}\left( \frac{\zeta _2 (t)}{ \zeta _1 (t)} \right) - \mathcal{M}_{+}}
	= \abs{ \CAL{M} \left(  \frac{a_{2} \zeta _2 (t)}{a_{2} \zeta _1 (t) -a_{1} \zeta _2 (t)}  \right) -1}
	\lesssim |x|^{ 2 \theta} \abs{\frac{\zeta _1(t)}{\zeta _2 (t)} - \frac{a_{1}}{a_{2}}}^{\theta}
	\label{le33:1}
\end{align}
for any $ \theta \in [0, 1]$.
By choosing suitable $\sigma_0$ and $r_1$, we can easily check that \eqref{exsi:1} 
satisfies Assumption \ref{A1}.
One also frequently utilizes $ \abs{ \zeta_{2}(t)} \lesssim |t|^{1-\lambda}$ which arises from \eqref{K13}.
The condition \eqref{K13-2} is only utilized in Lemma \ref{lem:hl1}.
Further, it is shown that there exists $\sigma(t)$ satisfying Assumption \ref{A1} in general (see \cite{GMT93, Na84, W69}).

Based on the Duhamel principle, we give the definition of solutions to \eqref{eq1} as follows:
\begin{Def}[$\mathcal{F}H^{\beta}$-solution]\label{def:sol}
Let $I \subset \R$ be an interval and fix $t_{0} \in I$. Set $\beta>0$.
We say a function $u(t,x): I \times \R^d \to \C$ is a $\mathcal{F}H^{\beta}$-solution to \eqref{eq1} on $I$ if $u \in C(I; L^2(\R^d))$
satisfies
\[
	u(t) = U(t, t_{0}) u(t_{0}) - i \int_{t_{0}}^{t} U(t, s) F(u(s))\, ds
\]
in $L^2(\R^d)$ for any $t \in I$ and $U(0, \cdot)u \in C(I; H^{0, \beta})$, where $H^{0, \beta}$ denotes the weighted $L^{2}$ space equipped with the norm $\norm{f}_{H^{0, \beta}} = \norm{(1 + |x|^{2})^{\beta/2} f}_{2}$.
\end{Def}

We now state the main theorems in this paper. Here we remark that the factorization \eqref{mdfm1} is only justified for the case $|t| \geq r_0$. Indeed, for the case of \eqref{exsi:1}, there may exists $t_0 \in [-r_0, r_0]$ such that $\zeta _j (t_0) = 0$. Hence we construct the modified scattering operator by the following three steps:
\begin{flushleft}
~~ 1. Construct the solution $u(t)$, $t\in (- \infty, -r_0]$ from $u_-$. \\ 
~~ 2. Construct the $u(0)$ from $u(-r_0)$ and $u(r_0) $ from $u(0)$ by the persistency argument. \\ 
~~ 3. Construct the $u_+$ and the solution $u(t)$, $t \in [r_0,\infty)$ from $u (r_0)$.   
\end{flushleft} 
The first main result is concerned with the modified wave operator as follows:

\begin{Thm}[Final state problem] \label{thm:1}
Assume $d/2 < \beta < \alpha< \min(1+p_c, d)$.
Then there exists $\ep_1>0$ such that the following assertion holds: For any $\varepsilon \in (0, \varepsilon_{1}]$ and 
for any $u_{-} \in H^{0, \alpha}$ with $\norm{u_{-}}_{H^{0, \alpha}} \leq \varepsilon$, 
there exists a unique $\mathcal{F}H^{\beta}$-solution to \eqref{eq1} on $(- \infty, -r_0]$.
Moreover the estimate
\[
	\norm{U(0, t)\left(u(t)- u_{p}(t)\right)}_{H^{0, \beta}} 
	\lesssim \varepsilon |t|^{-\mu}
\]
is valid for all $t \leq - r_0$, where $\mu>0$ is small enough and
\begin{align}
	u_p (t) = \mathcal{M} \left(  \frac{\zeta _2(t)}{\zeta _2 '(t)} \right) \mathcal{D} (\zeta _2 (t)) \widehat{u_{-}} \exp\(i \frac{\eta}{c_+} \left| \widehat{u_{-}}\right|^{p_{c}} \log |t|\), \quad c_+ = |a_{2}|^{\frac{1}{1-\lambda}}.
	\label{asym:1}
\end{align}
\end{Thm}

\begin{Rem}
If we replace the smallness assumption of $\norm{u_{-}}_{H^{0, \alpha}}$ by $\norm{\widehat{u_{-}}}_{ \infty} \leq \varepsilon$, then there exists $T \geq r_{0}$ such that \eqref{eq1} admits a unique $\mathcal{F}H^{\beta}$-solution on $(- \infty, -T]$ that behaves like $u_{p}$ as $t \rightarrow - \infty$ in the same topology. 
Although the replacement is the weaker condition, we need in the sequel the smallness in $H^{0, \alpha}$ to define the modified scattering operator. 
Hence we have already required the smallness.
\end{Rem}

\begin{Rem}
In \cite{Ka21} and \cite{KaMi23}, the time-weighted Strichartz estimate has been applied to discuss the final state problem. 
However, a use of the Strichartz estimate imposes the strong restriction of $ \lambda$.
We remove the restriction by considering a solution in the space associated with a generator of the Galilean transformation
$J(t) \coloneqq U(t, 0) x U(0, t)$.
\end{Rem}

\begin{Rem}
In Theorem \ref{thm:1}, the same assertion holds if we replace $u_{-}$ by $ \mathcal{M}_{+} u_{-}$ and 
$u_{p}(t)$ behaves like $U(t,0) \mathcal{F}^{-1} \exp\(-i \eta/c_{+} \left|\widehat{u_{-}}\right|^{p_{c}} \log |t|\) \widehat{u_{-}}$ as $t \rightarrow - \infty$ under suitable $u_{-}$.
\end{Rem}

We here define $B_{ \varepsilon}^{ \alpha} = \{ \phi \in H^{0, \alpha} \mid \norm{ \phi}_{H^{0, \alpha}} < \varepsilon \}$ for any $ \alpha$, $\varepsilon >0$.
Thanks to the persistence of regularity of solutions, see Lemma \ref{M-6:P1}, and the well-posedness in $L^{2}$ for \eqref{eq1}, we can extend the stopping time from $t = - r_{0}$ to $t= 0$ when the solution propagates.
Thus one can define the modified wave operator as follows:

\begin{Cor} \label{cor:waveop}
Under the assumption as in Theorem \ref{thm:1}, the modified wave operator $W_{-}$ on $B_{\varepsilon}$ is defined for any $ \varepsilon \in (0, \varepsilon_{1}]$.
Namely, 
for some $ c_{0} >0$, there exists a map
\begin{align*}
	\mathcal{W}_{-} \colon B^{\alpha}_{\varepsilon} \ni u_{-}  \mapsto u(0) \in B^{ \beta}_{c_{0} \varepsilon},
\end{align*}
where $ \varepsilon_{1}$ is given in Theorem \ref{thm:1}. 
\end{Cor}

Again using the persistence of regularity, for any $u(0) \in {B}_{c_0 \ep}^{\beta} $, there uniquely exists $\mathcal{F}H^{\beta}$-solution to \eqref{eq1} on $[0, r_0)$ such that for some $c_1 >0$ and any $t \in [0,r_0]$, $U(0,t)u(t) \in {B}_{c_1 \ep} ^{\beta}$ holds. 

Noting this fact, we state the second main result which is related to the existence of modified inverse wave operator.

\begin{Thm}[Modified scattering for the initial value problem] \label{thm:iv2}
Let $d/2 < \delta < \beta < \min(1+p_c, d)$.
Then there exists $\varepsilon_{2} >0$ such that the following assertion holds: For any $\varepsilon \in (0, \varepsilon_2]$ and $u_0 \in L^{2}$ satisfying $U(0, r_{0}) u_{0} \in H^{0, \beta}$ with $\norm{U(0, r_{0})u_{0}}_{H^{0, \beta}} \leq \varepsilon$, 
there exists a unique $\mathcal{F}H^{\beta}$-solution to \eqref{eq1} on $[r_{0}, \infty)$ under $u(r_{0})=u_0$.
Moreover, there exists a unique $u_{+} \in H^{0, \delta}$ such that
\begin{align}
	\norm{ \mathcal{F}^{-1} \mathcal{P}(t) \mathcal{F} \mathcal{M}_{+} U(0, t) u(t) - u_{+}}_{H^{0, \delta}} \lesssim{}& \varepsilon t^{- \mu}
	\label{thm:iv2a}
\end{align}
for any $t > r_{0}$, where 
\begin{align*}
	\mathcal{P}(t) = \exp \left( i \eta \int_{r_{0}}^{t} |\mathcal{F} \mathcal{M}_{+} U(0, \tau) u(\tau)|^{p_{c}} |\zeta_{2}(\tau)|^{- \frac{1}{1- \lambda}}\, d\tau  \right)
\end{align*}
if $d=1$, otherwise
\begin{align*}
	\mathcal{P}(t) = \exp \left( i \eta \int_{r_{0}}^{t} \( |\tau|^{- \chi} + |\mathcal{F} \mathcal{M}_{+} U(0, \tau) u(\tau)|^{2} \)^{\frac{p_{c}}{2}} |\zeta_{2}(\tau)|^{- \frac{1}{1- \lambda}}\, d\tau  \right)
\end{align*}
for some small $\chi>0$ and $\mu>0$. 
\end{Thm}

\begin{Rem}
The factor $ \mathcal{M}_{+}$ appears in the phase correction $\mathcal{P}(t)$, because we do not assume $a_{1} = 0$ unlike \cite{Ka21} and \cite{KM21}.
\end{Rem}

Similarly to the modified wave operator, 
we can construct the inverse wave operator as consequence of Theorem \ref{thm:iv2}.

\begin{Cor} \label{cor:inwa}
Under the assumption as in Theorem \ref{thm:iv2}, 
there exists $ \widetilde{\varepsilon}_{2} > 0$ such that for any $ \varepsilon \in (0, \widetilde{\varepsilon}_{2}]$,
the modified inverse wave operator $\CAL{W}_{+}^{-1}$ on $B_{ \varepsilon}^{ \beta}$ can be defined.
Namely, 
there exists a map
\begin{align*}
	\mathcal{W}_{+}^{-1} \colon B^{ \beta}_{ \varepsilon} \ni u(0)  \mapsto u_{+} \in H^{0, \gamma}.
\end{align*}
\end{Cor}

Combining Corollary \ref{cor:waveop} with Corollary \ref{cor:inwa}, the following holds: 
\begin{Cor} \label{cor:sca}
Let $d/2 < \delta < \alpha < \min \( d,1+p_c \)$ 
Then there exists $\varepsilon_{0} >0$ such that for any $ \varepsilon \in (0, \varepsilon_{0}]$, 
the modified scattering operator $\CAL{S}_+ \coloneqq \mathcal{W}_{+}^{-1} \mathcal{W}_{-}$  on $B_{ \varepsilon}^{ \alpha}$ can be defined. 
Namely, 
there exists a map
\begin{align*}
\CAL{S}_+ \colon B^{\alpha}_{ \varepsilon} \ni u_{-} \mapsto u_{+} \in H^{0, \delta}.
\end{align*} 
\end{Cor}

\subsection*{Notations}
We introduce some notations. 
For any $p \geq 1$, $L^p = L^{p}(\R^d)$ denotes the usual Lebesgue space on $\R^d$ equipped with the norm  $\norm{\, \cdot\, }_{p} \coloneqq \norm{\, \cdot\,}_{L^{p}}$.
Set $\J{a}=(1+|a|^2)^{1/2}$ for $a \in \C$ or $\R^d$. 
$\SCR{F}[u] = \widehat{u}$ is the usual Fourier transform of a function $u$ on $\R^d$ and $\SCR{F}^{-1}[u] = \check{u}$ is its inverse.  
Let $m$, $s \in \R$. 
The weighted Sobolev space and the homogeneous Sobolev space
on $\R^{d}$ are defined by $H^{m,s} = H^{m,s}(\R^d) = \{u \in \mathscr{S}' \mid \J{x}^s \J{i \nabla}^m u \in L^2 \}$ and 
$\dot{H}^s = \dot{H}^{s}(\R^d) = \{u \in \mathscr{S}' \mid |\nabla|^s u \in L^2\}$, respectively. 
Here $\mathscr{S}'$ is the space of tempered distributions.
We simply write $H^{m} = H^{m,0}$.
One denotes the smallest integer $n_0$ such that $n_0 \geq \delta$ by $\lceil \delta \rceil$ for any $\delta \in \R$.
$A \lesssim B$ denotes $A \leq CB$ for some constants $C>0$.

\subsection{Strategy on the proof of main results}

Before stating the strategy on the proof of main results, let us introduce several notations associated with the factorization  \eqref{mdfm1}, which is the key tool of our analysis. 
We denote 
\begin{align*}
\CAL{M}\left(  \frac{\zeta _2(t)}{\zeta _1 (t)}\right) = \CAL{M}_{+} \CAL{M}_{2}(t), \quad \CAL{M}_{+}= \CAL{M}\left( \frac{a_{2}}{a_{1}} \right), \quad \CAL{M}_{2}(t) = \CAL{M} \left(  \frac{a_{2} \zeta _2 (t)}{a_{2} \zeta _1 (t) -a_{1} \zeta _2 (t)}  \right) .
\end{align*} 
Then by \eqref{K13} and \eqref{le33:1}, we notice 
\begin{align} \label{K2/4-1}
\left| 
\left( \CAL{M}_2(t) \right) ^{\pm 1} -1 
\right| \lesssim
\begin{cases}
 |x|^{2 \theta} t^{- \theta (1- 2 \lambda)}, & \lambda >0, \\ 
 |x|^{2 \theta} t^{- \theta}, & \lambda  \leq 0   .
 \end{cases}
\end{align}
Hereafter one also uses the notation 
\begin{align*}
\CAL{M}_1 (t) = \CAL{M} \left(  \frac{\zeta _2(t)}{\zeta _2 '(t)} \right), \quad 
\mathcal{D}_{1}(t) = \mathcal{D}(\zeta_{2}(t)).
\end{align*}
Then \eqref{mdfm1} is rewritten as 
\begin{align}
	U(t,0) = \CAL{M}_{1}(t) \CAL{D}_{1}(t) \mathcal{F} \CAL{M}_{+} \CAL{M}_{2}(t).
	\label{mdfm}
\end{align}

\subsubsection{Strategy on the proof of Theorem \ref{thm:1}}
The strategy on the proof of Theorem \ref{thm:1} is to solve an integral equation associated with \eqref{eq1} around a prescribed asymptotic profile $u_{p}$, where
\begin{align}
u_p (t) ={}& \CAL{M}_1(t) \CAL{D}_{1}(t) \widehat{w} (t), \quad
	\widehat{w} (t)  = \widehat{u_{-}} \exp \left( i \frac{\eta}{c_{+}} |\widehat{u_{-}}|^{p_c} \log |t| \right).
	\label{thm:1ap}
\end{align}
We denote a new unknown $v$ by $v \coloneqq u - u_{p}$. 
Then one obtains the equation
\begin{align*}
	i \partial_t v -H(t) v 
	={}& F(u_{p}+v) - F(u_{p}) - \{ \( i \partial_t  -H(t) \) u_{p} -  F(u_{p}) \},
\end{align*}
which implies
\begin{align}
	i \partial_t \( U(0,t) v \) 
	={}& U(0,t) \{ F(u_{p}+v) - F(u_{p}) \} - \{ i \partial_t \(U(0,t) u_{p}\)  - U(0,t) F(u_{p}) \}.
	\label{intee}
\end{align}
Since $u_{p}(t)$ is an expected asymptotic profile of $u(t)$, we may assume that $\norm{v(t)}_{2} \rightarrow 0$ as $t \rightarrow -\infty$.
Therefore, applying \eqref{mdfm}, it follows from the Duhamel principle that $v(t)$ solves the integral equation
\begin{align}
	v(t) = -i \int_{-\infty}^{t} U(t,s) (F(u_{p}(s) + v(s)) - F(u_{p}(s)))\, ds + \mathcal{E}(t), \label{inte:1}
\end{align}
where $\mathcal{E}(t) = \mathcal{A}(t) + \mathcal{E}_{r}(t)$ with
\begin{align}
	\begin{aligned}
\CAL{A} (t) ={}&  -i \int_{-\infty}^{t} U(t,0) \CAL{M}_{2}(s)^{-1} \CAL{M}_{+}^{-1} \mathcal{F}^{-1} 
	\left( \frac{c_+ |s|}{|\zeta_2 (s)|^{1/(1- \lambda)} } -1 \right) F(\widehat{w}(s))\, \frac{ds}{c_{+} |s|}, \\[5pt]
\CAL{E}_{\mathrm{r}} (t) ={}& R(t) \widehat{w}(t) +i \int_{-\infty}^{t} U(t,s) R(s) F(\widehat{w}(s)) \frac{ds}{c_{+} |s|}, \\
R(t) ={}& \CAL{M}_1 (t) \CAL{D}_{1}(t) \left( \mathcal{F}\CAL{M}_2 (t) \mathcal{F}^{-1} -1 \right).
	\end{aligned}
	\label{inte:1remi}
\end{align}
Note that we give the derivation of \eqref{inte:1} in Appendix \ref{app:B}.
Conversely, if we find a solution to \eqref{inte:1}, then $u \coloneqq v + w$ is the solution to \eqref{eq1}
, that is $u(t)$ solves 
\begin{align} \label{K2/1-1}
	u(t) = u_{p}(t) +  i \int_{-\infty}^{t} U(t,s) (F(u(s)) - F(u_{p}(s)))\, ds + \mathcal{E}(t).
\end{align}
Hence we attempt to find a unique solution to \eqref{inte:1}, based on the contraction mapping argument in the complete metric space
\begin{align}
	\begin{aligned}
	X_{b, R} ={}& \left\{ \phi \in C\left( (- \infty, -r_{0}] ; L^2({\bf R}^d) \right) \relmiddle| \left\|\phi \right\|_{X_{b}} \leq R \right\}, \\
	\norm{\phi}_{X_{b}} ={}& \sup_{t \leq - r_{0}} \( |t|^{\frac{\beta}{2}(1-2\lambda) +b} \norm{\phi(t)}_{L^2(\R^d)}  + |t|^{b} \norm{U(0,t) \phi(t)}_{H^{0, \beta}(\R^d)} \), \\
	d(\phi, \psi) ={}& \norm{\phi - \psi}_{L^{\infty}(- \infty, -r_{0} ; L^2(\R^{d}))}
	\end{aligned}
	\label{K1}
\end{align}
for some $R>0$, where $b = b(\lambda, \beta, d)>0$ is small enough. 
If $ \norm{u_{-}}_{H^{0, \alpha}}$ is small enough and if the external term $ \mathcal{E}(t)$ satisfies $ \norm{ \mathcal{E}(t)}_{X_{b}} \leq  C_{0} \norm{u_{-}}_{\dot{H}^{0, \alpha}}$ for some $C_{0}>0$, then 
it enables us to show that \eqref{inte:1} admits a unique solution $u \in X_{b, R}$ with $R = 2C_{0} \norm{u_{-}}_{H^{0, \alpha}}$.
To this end, we need to estimate $F(u_{p}+v) - F(u_{p})$ in $X_{b, R}$ to show $\Psi (X_{b, R}) \subset X_{b, R}$, where $ \Psi$ is a map associated with \eqref{inte:1}.
Further, when $d=2$, $3$, the smallness of the exponent $p_{c}$ causes the lack of differentiability of the nonlinearity.
This means that if $X_{b, R}$ is equipped with the distance $d(\phi, \psi) = \norm{\phi - \psi}_{X_{b}}$, then it is difficult to show that $\Psi$ is the contraction map in $X_{b, R}$.
To overcome this difficultly, we prove the completeness of $X_{b, R}$ equipped with the distance in \eqref{K1} consisting of the weaker norm than $X_{b}$-norm.
This kind of technique is originally introduced by Kato \cite{TK87} (cf. \cite[p94]{C03}).

\subsubsection{Strategy on the proof of Theorem \ref{thm:iv2}} \label{ssec:1}

The proof of Theorem \ref{thm:iv2} consists of two parts.
The first part is to show that \eqref{eq1} admits a unique global $\mathcal{F}H^{\beta}$-solution that decays in time of the order $- d(1- \lambda)/2$ in $L^{\infty}$ (Theorem \ref{iv:thm1} below), and the second part is to specify the final data $u_{+}$ in $H^{0, \delta}$. 
As mentioned above, under $u(0) = u_{0} \in H^{ \beta} \cap H^{0, \beta}$ with $\beta > \max(1, d/2)$,  the existence of such a solution and $u_{+}$ in $L^{2}$ are shown in \cite{KM21}.
To explain the strategy on their proof of the first part, we focus attention on the following inequality given by \eqref{mdfm}: 
\begin{align*}
	\norm{u(t)}_{\infty} \lesssim |t|^{-\frac{d}{2}(1 -\lambda)} \norm{\mathcal{F} \mathcal{M}_{+} U(0, t) u(t)}_{\infty}
	+ o(|t|^{-\frac{d}{2}(1- \lambda)}) 
\end{align*}
for any $t \geq r_{0}$.
According to the above inequality, we have to show the boundedness of $\norm{\mathcal{F} \mathcal{M}_{+} U(0, t) u(t)}_{\infty}$ to get the desired solution.
This claim can be proven by constructing the solution $u$ that $|J(t)|^{\beta} u(t)$ increase in time of the small order depending on the size of $u_{0}$.
This kind of technique is originally developed by \cite{HN98}.
In \cite{KM21}, $u_{0} \in H^{ \beta}$ is required to control $ \norm{u(t)}_{ \infty}$ near $t = 0$, since they construct the solution $u(t)$ on $[0, \infty)$. 
To avoid a use of the condition, we construct the solution on $[r_{0}, \infty)$.
They also employ \eqref{mdfm1} instead of \eqref{mdfm}, because the factor $\mathcal{M}_{+}$ does not appear due to the assumption $a_{1} = 0$ in \cite{KM21}.
Hence, for the reader's convenience, we give a sketch of the proof of Theorem \ref{iv:thm1} in Appendix \ref{app:A}, although it has been essentially proven in \cite{KM21}.

The novelty of the paper is to find the final data $u_{+}$ in $H^{0, \delta}$ with $d/2<\delta< \beta$.
To state our technique, we only handle $d=2$, $3$.
Following \cite{HN98}, from \eqref{eq1}, we derive
\begin{align*}
	i\partial_t (\mathcal{F} \mathcal{M}_{+} U(0, t)u) 
	={}& \eta |\zeta_{2}(t)|^{-\frac{1}{1- \lambda}} |\mathcal{F} \mathcal{M}_{+} U(0, t)u|^{p_{c}}\mathcal{F} \mathcal{M}_{+} U(0, t)u
	+ o(t^{-1}).
\end{align*}
Roughly speaking, the first term in the above right-hand side is $\mathcal{O}(t^{-1})$, that is, it is non-integrable with respect to time on $[r_{0}, \infty)$.
Hence we introduce a phase modification
\begin{align*}
	 v_{p} \coloneqq{}&  \mathcal{F}^{-1} \mathcal{P}(t) \mathcal{F} \mathcal{M}_{+} U(0, t)u, \\
	 \mathcal{P}(t) \coloneqq{}& \exp \left( i \eta \int_{r_{0}}^{t} \( |\tau|^{- \chi} + |\mathcal{F} \mathcal{M}_{+} U(0, \tau) u(\tau)|^{2} \)^{\frac{p_{c}}{2}} |\zeta_{2}(\tau)|^{- \frac{1}{1- \lambda}}\, d\tau  \right)
\end{align*}
for some $\chi>0$.
By the modification, we obtain
\begin{align*}
	i \partial_t (\mathcal{F}v_{p}) ={}& \eta |\zeta_{2}(t)|^{-\frac{1}{1- \lambda}} 
	\mathcal{P}(t) \( |\mathcal{F} \mathcal{M}_{+} U(0, t)u(t)|^{p_{c}} - \( |t|^{- \chi} + |\mathcal{F} \mathcal{M}_{+} U(0, t)u(t)|^{2} \)^{\frac{p_{c}}{2}} \) \\
	&{} \times \mathcal{F} \mathcal{M}_{+} U(0, t)u(t) + o(t^{-1}).
\end{align*}
This implies that there exists $u_{+} \in H^{0, \delta}$ such that $v_{p} \rightarrow u_{+}$ in $H^{0, \delta}$ as $t \rightarrow \infty$.
Owing to the lack of differentiability of the nonlinearity, the argument in \cite{HN98} is not directly applicable, since we need to find the final data in $H^{0, \delta}$.
To overcome the difficulty, we take the above phase modification with the factor $|\tau|^{- \chi}$, which is firstly used by Ginibre-Ozawa \cite{GO93}.
The key of the proof is to see that the first term of the above right-hand side, which comes from the modification, is negligible for large time.
To deal with this issue, we develop fundamental inequalities necessary to extract the time decay of the first term. 

\subsection{Organization of this paper}

The rest of the paper is organized as follows:
In Section \ref{Sec:2}, we collect several ingredients to construct the modified scattering operator.
One first states properties for the Galilean operator $J(t)$ and estimates for the fractional calculus. 
We also prove a persistent property for the regularity of solutions necessary to connect the solution near $t = 0$.
Nonlinear estimates involving fractional derivatives are discussed at the end of this section.
We next give the proof of Theorem \ref{thm:1} and construct the modified wave operator in Section \ref{Sec:4}.
In Section \ref{Sec:5}, key inequalities for modified phase corrections to get rid of the lack of differentiability of the nonlinearity are developed, and the proof of Theorem \ref{thm:iv2} will be given.
Moreover, we shall construct the modified inverse wave operator and finally define the modified scattering operator at the end of this section. 
Appendix \ref{app:A} provides us to prove Theorem \ref{iv:thm1}, which is the former part of the proof of Theorem \ref{thm:iv2}. 
The derivation of the integral equation \eqref{inte:1} is also treated in Appendix \ref{app:B}.
For the sake of simplicity, we give all of statements in the both case $p_{c} \geq 1$ and $p_{c} <1$, but usually deal with $p_{c} <1$ only because the other case are easier or similar.
Remark that $p_c \geq 1$ holds when $d=1$, $2$ or $d=3$ and $\lambda \in [1/3, 1/2)$ and also if $d=3$ and $\lambda \in (-1/3, 1/3)$, then $p_c < 1$. 
Because $1+ p_c \leq d/2$ holds in the case $d=3$ and $\lambda \leq -1/3$, we can not rely on the approach in the use of the Sobolev inequality and hence we omit to discuss this case. Furthermore, the case $$ \lambda \geq 0 $$ is only treated afterward since the case $\lambda<0$ is same as $\lambda = 0$.


\section{Preliminaries} \label{Sec:2}

\subsection{Some tools}

We summarize several tools frequently used throughout this paper.
Let us first introduce the operator $J(t) = U(t, 0) x U(0, t)$ which is often called the Galilean operator.
Besides, define $|J(t)|^{ \gamma} = U(t, 0) |x|^{ \gamma} U(0, t)$ for any $ \gamma \geq 0$.
The operator plays a important role in this paper. 
Noting that 
\begin{align*}
|J(t)|^{\gamma} ={}& \CAL{M}_1(t) \CAL{D}_1(t) \CAL{F} \CAL{M}_+ \CAL{M}_2(t) |x| ^{\gamma} \CAL{M}_2(t) ^{-1} \CAL{M}_+^{-1} \CAL{F}^{-1}  \CAL{D}_1(t)^{-1} \CAL{M}_1(t)^{-1}  , \\ 
|J(t)|^{\gamma} R(t) 
	={}& U(t,0) |x|^{\gamma} U(0,t) \cdot \CAL{M}_1(t) \CAL{D}_1(t) \left( \CAL{F}\CAL{M}_2(t) \CAL{F}^{-1} -1   \right)  ,
\end{align*} 
we have the following properties of the operator: 
\begin{Lem} \label{lem:k1}
For all $t \in {\bf R}$ and $\gamma>0$, it holds that
\begin{align}
	|J(t)|^{\gamma} ={}& 
	 \left|\zeta_2(t)\right|^\gamma \mathcal{M}_{1}(t) | \nabla|^\gamma \mathcal{M}_{1}^{-1}(t), \label{opa:1}\\
	|J(t)|^{\gamma} R(t)  
	={}& R(t) | \nabla|^{ \gamma}. \label{opa:3}
\end{align}
\end{Lem}

We often employ the following well-known estimates:
\begin{Lem}[e.g. {\cite[Lemma 4.6]{KaMi23}}] \label{lem:nl1}
Let $\alpha>0$ and $n \in \Z$. If $\alpha \geq 1$, then
\begin{align}
	\left| |z|^{\alpha-n}z^n - |w|^{\alpha-n}w^n \right| \leq C |n| \left(|z|^{\alpha-1} + |w|^{\alpha-1}\right) |z-w| \label{nl:1}
\end{align}
for any $z$, $w \in \C$. Further, when $ \alpha \in (0,1)$, the following holds:
\begin{align}
	\left| |z|^{\alpha-n}z^n - |w|^{\alpha-n}w^n \right| \leq C |z-w|^{\alpha}. \label{nl:2}
\end{align}
\end{Lem}

The subsequent estimates are the fractional Leibniz rule and chain rule.

\begin{Lem}[\cite{KaP88, GraO14}]
\label{lem_lei}
Let $s>0$, $1<r<\infty$, $1<p_1,p_2,q_1,q_2\leq \infty$ and $1/r=1/p_1+1/q_1=1/p_2+1/q_2$. Then we have the following fractional Leibniz rule: 
$$
\norm{|\nabla|^s(fg)}_{r} \lesssim \( \norm{|\nabla|^sf}_{p_1}\norm{g}_{q_1}+\norm{f}_{p_2}\norm{|\nabla|^sg}_{q_2}\).
$$
\end{Lem}

\begin{Lem}[\cite{Vi07}]
\label{lem_hcha}
Let $F$ be a H\"older continuous function of order  $0< \rho <1$. Suppose that $0<\sigma<\rho$, $1<p<\infty$ and $\sigma/\rho<s<1$. Then
\begin{align*}
\norm{|\nabla|^\sigma F(u)}_{L^p} \lesssim \norm{|u|^{\rho-\sigma/s}}_{L^{p_1}}\norm{|\nabla|^su}_{L^{p_2\sigma /s}}^{\sigma/s}
\end{align*}
provided $1/p = 1/p_1 + 1/p_2$ and $(1 - \sigma/(\rho s))p_1 >1$.
\end{Lem}

\subsection{Persistence of regularity of solutions}
In Section \ref{Sec:4}, we construct the $\CAL{F}H^{\beta}$-solution $u(t)$, $t \in (- \infty, -r_0]$ from $u_-$, while the construction scheme strongly rely on \eqref{mdfm1} which is the factorization of $U(t,s)$. In order to construct the modified wave operator, we need to construct the nonlinear propagator  $u(-r_0) \mapsto  u(0)$ and in $[-r_0, 0]$, we can not rely on the approach due to factorization. 
Similarly, in order to define the inverse wave operator in Section \ref{Sec:5}, one needs the nonlinear propagator $u(0) \mapsto u(r_0)$. 
In this subsection, these propagators are established using the persistence of regularity.
We here show the persistence result in a generalized form, but, in our case, it is enough to consider the case where $I = [-r_0 , 0]$, $t_0 = -r_0$ and $u_0 = u(-r_0)$, ({resp.} $I=[0,r_0]$, $t_0 =0$ and $u_0 = u(0)$) in Proposition \ref{pro:ext1} as follows: 
\begin{Prop} \label{pro:ext1}
Let $I$ be a bounded interval and fix $t_{0} \in I$ and $(Q_0, R_0) = ( 4(2+p_c)/(dp_c) , 2+p_c)$.
Assume $U(0, t_{0}) u_{0} \in H^{0, \beta}$ with $\| U(0, t_{0}) u_{0}  \|_{H^{0,\beta}} $ small enough compared to $|I|$.
Then \eqref{eq1} has a unique solution $u \in C(I; L^{2}) \cap L^{Q_0} (I; L^{R_0}), |J (\cdot) |^{\beta} u \in C(I; L^{2}) \cap L^{Q_{0}}(I; L^{R_{0}})$ with
\begin{align}
	\norm{u}_{L^{ \infty}(I; L^{2}) \cap L^{Q_{0}}(I; L^{R_{0}})} 
	+ \norm{|J( \cdot)|^{ \beta} u}_{L^{ \infty}(I; L^{2}) \cap L^{Q_{0}}(I; L^{R_{0}})} 
	\leq {}& A_0 \norm{U(0, t_{0}) u_{0}}_{H^{0, \beta}}, \label{pere:1}
\end{align}
where $A_{0}>0$ depends on $d$, $ \lambda$, $I$ and $ \norm{u_{0}}_{2}$ and monotonically increases for $\norm{u_{0}}_{2}$.
\end{Prop}

\begin{Rem}
We only employ Proposition \ref{pro:ext1} to construct the propagators on $[-r_{0}, r_{0}]$.
\end{Rem}

Let us divide several steps into the construction of the nonlinear propagators.
Following the pioneer work of Tsutsumi \cite{YT87},
Carles \cite{Ca11} proves the global well-posedness in $L^{2}$ by the local-in time Strichartz estimate and the conservation of $\norm{u(t)}_{2}$. 

\begin{Lem}[{cf. \cite[Proposition 1.5]{Ca11}}] \label{K-5/14-L1}
Let $I$ be a bounded interval and fix $t_{0} \in I$ and $(Q_0, R_0) = ( 4(2+p_c)/(dp_c) , 2+p_c)$.
Assume $u_{0} \in L^2$.
Then \eqref{eq1} has a unique global solution $u \in C(I; L^{2}) \cap L^{Q_0} (I; L^{R_0}) $ with $u(t_{0}) = u_{0}$ such that 
\begin{align}
	\norm{u}_{L^{ \infty}(I; L^{2}) \cap L^{Q_{0}}(I; L^{R_{0}})} 
	\leq{}& A_1 \norm{u_{0}}_{2}, 
	\label{yt87:1}
\end{align}
where $A_{1}>0$ depends on $d$, $ \lambda$, $I$ and $\norm{u_{0}}_{2}$ and monotonically increases for $\norm{u_{0}}_{2}$.
\end{Lem}
By means of the nonlinear estimates proven by the fractional calculus, 
we can obtain the following persistent property for the regularity of solutions. 
We here recall $|J(t)|^{\beta} = U(t,0) |x|^{ \beta} U(0,t)$ in $L^{2}$.
Remark that a similar property has been discussed in \cite{NP09}.

\begin{Lem}[Persistence of regularity of solutions] \label{M-6:P1}
Let $u$ be the solution as in Lemma \ref{K-5/14-L1} with the same interval $I$ under $u(t_{0}) = u_{0}$ and suppose $U(0, t_{0})u_{0} \in H^{0, \beta}$.
Then $|J (\cdot) |^{\beta} u \in C(I; L^2) \cap L^{Q_{0}}(I; L^{R_{0}})$ and 
it holds that
\begin{align}
	\norm{|J( \cdot)|^{ \beta} u}_{L^{ \infty}(I; L^{2}) \cap L^{Q_{0}}(I; L^{R_{0}})} 
	\leq{}&  A_2 \norm{U(0, t_{0}) u_{0}}_{H^{0, \beta}} \label{pere:2}, 
\end{align}
where $A_{2}>0$ depends on $d$, $ \lambda$, $I$ and $ \norm{u_{0}}_{2}$ and monotonically increases for $\norm{u_{0}}_{2}$.
\end{Lem}

In order to prove Lemma \ref{M-6:P1}, we need the following nonlinear estimate:


\begin{Lem} \label{M-6:L1}
Let $(Q_0, R_0) = ( 4(2+p_c)/(dp_c) , 2+p_c )$ as in Lemma \ref{K-5/14-L1}.
Then the estimate
\begin{align*}
	\norm{|J(\cdot)|^{\beta}|u|^{p_c}u}_{L^{Q_0'}(I; L^{R_0'})}
	\lesssim{}& |I|^{\theta} \norm{u}_{L^{Q_0}(I; L^{R_0})}^{p_c} \norm{|J(t)|^{\beta}u}_{L^{Q_0}(I; L^{R_0})} 
\end{align*}
holds, where $\theta = 1 - d p_c/4 = 1 - 1/2(1- \lambda) \geq 0$. 
\end{Lem}

\begin{proof}[\textit{\bf Proof of Lemma \ref{M-6:L1}}]
We only consider the case $p_{c} <1$ since the other cases are similar by a use of the fractional chain rule (e.g. \cite{CW91, K95}).
Set $\phi = \mathcal{M}_{1}^{-1}(t) u$.
By \eqref{opa:1}, we have
\begin{align*}
\left\| 
|J(t)|^{\beta} |u(t)|^{p_c} u(t)
\right\|_{R_0} = |\zeta _2(t)|^{\beta} \left\| 
|\nabla|^{\beta } |\phi (t)|^{p_c} \phi (t)
\right\|_{R_0}. 
\end{align*} 
Straightforward calculation gives
\begin{align*}
	\nabla \( |\phi|^{p_c} \phi \) 
	={}& \(\frac{p_c}{2} + 1\) |\phi|^{p_c} \nabla \phi  
	+ \frac{p_c}{2}  |\phi|^{p_c-2} \phi^2 \overline{\nabla \phi} 
	\eqqcolon B_{1}( \phi) \nabla \phi + B_{2}( \phi) \overline{\nabla \phi}. 
\end{align*}
By Lemma \ref{lem_lei}, we obtain
\begin{align*}
	\norm{| \nabla|^{\beta -1} B_1(\phi) \nabla \phi}_{R_0'} 
	\lesssim{}& \norm{|\nabla|^{\beta -1}B_1(\phi)}_{R_0/p_c} \norm{\nabla \phi}_{R_0} 
	+ \norm{\phi}_{R_0}^{p_c} \norm{|\nabla|^{\beta} \phi}_{R_0}.
\end{align*}
Note that $1/R_0' = p_c/R_0 + 1/R_0$.
Applying Lemma \ref{lem_hcha}, one has
\begin{align*}
	\norm{|\nabla|^{\beta -1}B_1(\phi)}_{R_0/p_c}
	\lesssim{}& \norm{|\phi|^{p_c - \frac{\beta-1}{s}}}_{R_1} \norm{|\nabla|^{s} \phi}_{\frac{\beta-1}{s} R_2}^{\frac{\beta-1}{s}}
\end{align*}
for any $s \in \( (\beta-1)/p_c, 1\)$, where 
\[
	p_c/R_0 = 1/R_1 + 1/R_2, \quad 1/R_1 = (p_c-(\beta-1)/s)/R_0, \quad 1/R_2 = (\beta-1)/s \cdot 1/R_0.
\]
Hence, we see from the Gagliardo-Nirenberg inequality (see \cite[p168]{BL76}) that 
\begin{align*}
	\norm{\nabla \phi}_{R_0} \lesssim \norm{ \phi}_{R_0}^{1-\theta_1} \norm{|\nabla|^{\beta} \phi}_{R_0}^{\theta_1}, \quad
	\norm{|\nabla|^{s} \phi}_{R_{0}} \lesssim \norm{\phi}_{R_{0}}^{1-\theta_2} \norm{|\nabla|^{\beta} \phi}_{R_{0}}^{\theta_2},
\end{align*}
where $\theta_1 = 1/\beta \in (0,1)$ and $ \theta_2 = s/\beta \in (0,1)$.
Collecting the above estimates, 
\begin{align*}
	\norm{|\nabla|^{\beta -1}B_1(\phi)}_{p_c/R_0} \norm{\nabla \phi}_{R_{0}} 
	\lesssim{}& \norm{\phi}_{R_{0}}^{p_c} \norm{|\nabla|^{\beta} \phi}_{R_{0}}.
\end{align*}
holds. Since the estimate $\norm{B_2(\phi) \overline{\nabla \phi}}_{R_{0}}$ is the same, we obtain
\begin{align*}
	\norm{|J(t)|^{\beta}|u(t)|^{p_c}u(t)}_{R_0'}
	\lesssim \norm{u(t)}_{R_{0}}^{p_c} \norm{|J(t)|^{\beta}u(t)}_{R_{0}}.
\end{align*}
By the H\"older inequality, one has
\begin{align*}
	\norm{|J(\cdot)|^{\beta}|u|^{p_c}u}_{L^{Q_0'}(I, L^{R_0'})}
	\lesssim{}& |I|^{\theta} \norm{u}_{L^{Q_0}(I; L^{R_0})}^{p_c} \norm{|J(t)|^{\beta}u}_{L^{Q_0}(I; L^{R_0})}, 
\end{align*}
where $\theta = 1/Q_0' - (p_c+1)/Q_0 = 1- d p_c/4$, as desired.
\end{proof}

\begin{proof}[\textit{\bf Proof of Lemma \ref{M-6:P1}}]
Let us only consider the forward direction in time, because the backward direction is same.
Set $T_{0} = \sup I - t_{0}$.
From the standard well-posedness theory, there exists a unique local $\mathcal{F}H^{ \beta}$-solution $u$, $J(\cdot) u \in C([t_{0}, t_{0}+T]; L^{2}) \cap L^{Q_{0}}(t_{0}, t_{0}+T; L^{R_{0}})$ to \eqref{eq1} with $u(t_{0}) = u_{0}$ for some $T \leq T_{0}$.
We here denote the maximal existence time of its solution by $T_{{\rm max}}$.
Let us prove $T_{\max}= T_{0}$ by a contradiction. 
Suppose that $T_{\max} < T_{0}$.
Take arbitrarily small $\varepsilon$, $ \delta >0$ with $\delta > \varepsilon$. We set $T_{\tau} = T_{\rm{max}} - \tau < T_{\rm{max}}$ and $I_{\tau} = [t_{0}, t_{0} + T_{ \tau}]$ for any small $\tau >0$. 
By means of the local-in-time Strichartz estimate,
we see from Lemma \ref{M-6:L1} and \eqref{yt87:1} that
\begin{align*}
	&{}\norm{|J( \cdot)|^{ \beta} u}_{L^{ \infty}(I_{ \varepsilon}; L^{2})} + \norm{|J( \cdot)|^{ \beta} u}_{L^{Q_{0}}(I_{ \varepsilon}; L^{R_{0}})} \\
	\leq{}& 2C_{0} \norm{|J(t_{0})|^{ \beta}u_{0}}_{2} \\
	&{}+ C \norm{u_{0}}_{2}^{p_c} 
	\( |T_{ \delta}|^{ \theta} \norm{|J( \cdot)|^{\beta}u}_{L^{Q_0}(I_{ \delta}; L^{R_0})} + | \delta - \varepsilon|^{\theta} \norm{|J( \cdot)|^{\beta}u}_{L^{Q_0}(t_{0} +T_{ \delta}, t_{0}+ T_{\varepsilon}; L^{R_0})}\) \\
	\leq{}& C_{\delta}(u_{0}) + a_{1} \delta^{\theta} \norm{u_{0}}_{2}^{p_c} \norm{|J( \cdot)|^{\beta}u}_{L^{Q_0}(I_{ \varepsilon}; L^{R_0})}, 
\end{align*}
where $\theta = 1- d p_c/4$.
Hence taking $ \delta_{0} >0$ such that $a_{1} | \delta_{0}|^{\theta} \norm{u_{0}}_{2}^{p_c} \leq 1/2$, one has
\begin{align*}
	\norm{|J( \cdot)|^{ \beta} u}_{L^{ \infty}(I_{ \varepsilon}; L^{2})} + \norm{|J( \cdot)|^{ \beta} u}_{L^{Q_{0}}(I_{ \varepsilon}; L^{R_{0}})} \leq C_{ \delta_{0}}(u_{0}).
\end{align*}
Letting $ \varepsilon \downarrow 0$, we have the contradiction and thus $T_{\max} = T_{0}$.
This yields $|J( \cdot)|^{ \beta}u \in C( I; L^{2}) \cap L^{Q_{0}}(I; L^{R_{0}})$.

As for the estimate \eqref{pere:2}, similarly to the above, there exists $ \tau_{0}>0$ such that
\begin{align*}
	\norm{|J( \cdot)|^{ \beta} u}_{L^{ \infty}(I'_{\tau_{0}}; L^{2})} +  \norm{|J( \cdot)|^{ \beta} u}_{L^{Q_{0}}(I'_{ \tau_{0}}; L^{R_{0}})} \leq 2C_{0} \norm{|J(t_{0})|^{ \beta}u_{0}}_{2},
\end{align*}
where $I'_{ \tau_{0}} = [t_{0}, t_{0} + \tau_{0}] \subset I$.
Since $ \tau_{0}$ only depends on 
$\norm{u_{0}}_{2}$, 
thanks to the conservation of $\norm{u(t)}_{2}$,
we can iterate the above argument with the same $\tau_{0}$ and hence the desired assertion holds. 
This completes the proof.
\end{proof}
\begin{proof}[\textit{\bf Proof of Proposition \ref{pro:ext1}}]
It is the consequence of Lemma \ref{K-5/14-L1} and Lemma \ref{M-6:P1}. 
\end{proof}


\subsection{Nonlinear estimates} \label{Sec:3}

In this section, we shall deal with the following estimates with regard to the difference of the nonlinearity:

\begin{Lem} \label{lem:34}
Let $d/2 < \beta < \min\(d, 1+p_c\)$. Then the followings hold:
If $d=1$, then
\begin{align*}
	\norm{F(\phi) - F(\psi)}_{\dot{H}^{\beta}}  
	\lesssim{}& \left(\norm{\phi - \psi}_{\infty} + \norm{\psi}_{\infty} \right)^{p_c} \norm{\phi-\psi}_{\dot{H}^\beta} \\
	&{}+ \left(\norm{\phi - \psi}_{\infty} + \norm{\psi}_{\infty} \right)^{p_c-1} \norm{\phi-\psi}_{\infty}\norm{\psi}_{\dot{H}^\beta}.
\end{align*}
When $d=2$, $3$, the estimate
\begin{align*}
	&{}\norm{F(\phi) - F(\psi)}_{\dot{H}^\beta} \\
	\lesssim{}& \norm{\phi}_{\infty}^{p_c+1-\beta} \norm{\nabla \phi}_{2 \beta}^{\beta-1} \norm{\nabla (\phi - \psi)}_{2 \beta} 
	+ \norm{\phi}_{\infty}^{p_c}  \norm{\phi - \psi}_{\dot{H}^{\beta}} 
	+ \norm{\phi - \psi}_{\infty}^{p_{c}} \norm{ \psi}_{\dot{H}^{ \beta}} \\
	&{}+ \norm{\phi - \psi}_{\infty}^{p_{c} - \gamma( \beta -1)} 
	\( \norm{\phi}_{\infty}^{\gamma -1} \norm{\nabla \phi}_{2\beta}
	+ \norm{\psi}_{\infty}^{\gamma -1} \norm{\nabla \psi}_{2\beta} \)^{ \beta -1} \norm{\nabla \psi}_{2 \beta}
\end{align*}
holds for any $\gamma \in (1,  p_{c}/ (\beta-1))$ if $p_{c}<1$, otherwise
\begin{align*}
	\norm{F( \phi) - F( \psi)}_{\dot{H}^{\beta}} 
	\lesssim{}& \norm{\phi}_{\infty}^{p_c+1-\beta} \norm{\nabla \phi}_{2\beta}^{\beta-1} \norm{\nabla (\phi - \psi)}_{2\beta} 
	+ \norm{\phi}_{\infty}^{p_c}  \norm{\phi - \psi}_{\dot{H}^{\beta}} \\
	&{}+ \( \norm{\phi}_{\infty} + \norm{\psi}_{\infty} \)^{p_{c} -1} \norm{\phi - \psi}_{\infty} \norm{ \psi}_{\dot{H}^{ \beta}} \\
	&{}+ \( \norm{\phi}_{\infty} + \norm{\psi}_{\infty} \)^{ (p_{c} -1)(2 - \beta)} \norm{\phi - \psi}_{\infty}^{2- \beta} \\
	&{}\times \( \norm{ \phi}_{\infty}^{p_{c} -1} \norm{ \nabla \phi}_{2\beta} 
	+ \norm{ \psi}_{\infty}^{p_{c} -1} \norm{ \nabla \psi}_{2\beta} \)^{ \beta-1} \norm{\nabla \psi}_{2 \beta}.
\end{align*}
\end{Lem}
\begin{Rem}
Connecting Lemma \ref{lem:34} with $\psi=0$ and the Gagliardo-Nirenberg inequality, we easily verify $\norm{F(\phi)}_{\dot{H}^{\beta}} \lesssim \norm{\phi}_{\infty}^{p_c} \norm{\phi}_{\dot{H}^\beta}$.
\end{Rem}

\begin{proof}
The estimate in $d=1$ is immediate from a use of the representation of $\dot{H}^{\beta}$-norm via the Gagliardo semi-norm (see \cite[Proposition 3.4]{DPV12}), after the calculation as follows:
\begin{align*}
	F(\phi) - F(\psi) 
	={}& (\phi - \psi) \int_0^1\diffp{F}{z}(\psi + \theta (\phi - \psi))\, d \theta 
	+ (\overline{\phi - \psi}) \int_0^1\diffp{F}{\overline{z}}(\psi + \theta (\phi - \psi))\, d \theta.
\end{align*}
Let us consider the case $d=2$, $3$. 
Remark that
\begin{align}
\begin{aligned}
	\nabla \( |\phi|^{p_c} \phi - |\psi|^{p_c} \psi \) 
	={}& \(\frac{p_c}{2} + 1\) \( |\phi|^{p_c} \nabla (\phi -\psi) + \( |\phi|^{p_c} - |\psi|^{p_c} \) \nabla \psi \) \\
	&{}+ \frac{p_c}{2} \( |\phi|^{p_c-2} \phi^2 \nabla (\overline{\phi -\psi}) + \( |\phi|^{p_c-2} \phi^2 - |\psi|^{p_c-2} \psi^2 \) \overline{\nabla \psi} \).
\end{aligned}
	\label{nes32:1}
\end{align}
We here mimic the argument of the proof of \cite[Proposition 4.3]{KaMi23} to estimate the first term of the right-hand side of \eqref{nes32:1}.
Remark that it comes from \cite[Lemma 4.9]{KaMi23} that
\begin{align}
\begin{aligned}
	\norm{|\nabla|^{\beta-1} \left( |\varphi|^{\gamma-n}\varphi^n \right)}_{\frac{2\beta}{\beta-1}} 
	\lesssim{}& \norm{\varphi}^{\gamma + 1-\beta}_{\infty} \norm{\nabla \varphi}_{2\beta}^{\beta-1}
\end{aligned}
	\label{ino:1}
\end{align}
for any $n \in \Z$ and all $\gamma > \beta -1$.
Applying Lemma \ref{lem_lei} and \eqref{ino:1},
we reach to
\begin{align*}
	\norm{|\phi|^{p_c} \nabla (\phi -\psi)}_{\dot{H}^{\beta-1}} 
	\lesssim{}& \norm{|\nabla|^{\beta-1} |\phi|^{p_c}}_{\frac{2\beta}{\beta-1}} \norm{\nabla (\phi - \psi)}_{2\beta} 
	+ \norm{ |\phi|^{p_c}}_{\infty} \norm{|\nabla|^{\beta} (\phi - \psi)}_{2} \\
	\lesssim{}&  \norm{\phi}_{\infty}^{p_c+1-\beta} \norm{\nabla \phi}_{2\beta}^{\beta-1} \norm{\nabla (\phi - \psi)}_{2\beta} 
	 + \norm{\phi}_{\infty}^{p_c}  \norm{\phi - \psi}_{\dot{H}^{\beta}}.
\end{align*} 
Let us next go to the second term of \eqref{nes32:1}. 
we see from Lemma \ref{lem_lei} that 
\begin{align}
\begin{aligned}
	&{}\norm{\( |\phi|^{p_c} - |\psi|^{p_c} \) \nabla \psi}_{\dot{H}^{\beta-1}} \\
	\lesssim{}& \norm{ |\phi|^{p_c} - |\psi|^{p_c} }_{\infty} \norm{ \psi}_{\dot{H}^{ \beta}}
	+ \norm{ | \nabla|^{ \beta -1} \( |\phi|^{p_c} - |\psi|^{p_c} \) }_{\frac{2\beta}{\beta-1}} 
	\norm{\nabla \psi}_{2 \beta}.
\end{aligned}
	\label{es:823}
\end{align}
Hereafter, one shall treat the case $p_{c}<1$.
In view of the concavity of $|z|^{p_{c}}$, 
\begin{align*}
	\norm{|\phi|^{p_c} - |\psi|^{p_c}}_{\infty} \leq \norm{\phi - \psi}_{\infty}^{p_c}
\end{align*}
holds. 
Hence, the Gagliardo-Nirenberg inequality leads to 
\begin{align}
\begin{aligned}
	\norm{|\nabla|^{\beta-1} \left( |\phi|^{p_c} - |\psi|^{p_c} \right)}_{\frac{2\beta}{\beta-1}} 
	\lesssim{}& \norm{\phi - \psi}_{\infty}^{p_{c} \(1- \frac{\beta-1}{s_0}\)} \\
	&{}\times \( \norm{ |\nabla|^{s_0}  |\phi|^{p_c} }_{\frac{2\beta}{s_0}}
	+ \norm{ |\nabla|^{s_0} |\psi|^{p_c} }_{\frac{2\beta}{s_0}} \)^{\frac{\beta-1}{s_0}},
\end{aligned}
	\label{app:n2}
\end{align}
where $s_0 = p_{c} - \kappa$ for any $ \kappa \in (0, p_{c} - \beta+1)$.
By Lemma \ref{lem_hcha}, we deduce that 
\begin{align*}
	\norm{|\nabla|^{s_0} |\varphi|^{p_{c}} }_{\frac{2\beta}{s_0}} 
	\lesssim{}& \norm{\varphi}_{\infty}^{p_{c} - \frac{s_0}{s}} \norm{|\nabla|^s \varphi}_{\frac{2\beta}{s}}^{\frac{s_0}{s}}
\end{align*}
for any $s \in \(s_0/p_{c}, 1\)$.
We then see from the Gagliaro-Nirenberg inequality that 
\begin{align*}
	\norm{|\nabla|^{s_0} |\varphi|^{p_{c}}}_{\frac{2\beta}{s_0}} 
	\lesssim{}& \norm{\varphi}_{\infty}^{p_{c} -s_0} \norm{\nabla \varphi}_{2\beta}^{s_0}.
\end{align*}
Plugging the above into \eqref{app:n2}, one has
\begin{align*}
	\norm{|\nabla|^{\beta-1} \left( |\phi|^{p_c} - |\psi|^{p_c} \right)}_{\frac{2\beta}{\beta-1}} 
	\lesssim{}& \norm{\phi - \psi}_{\infty}^{p_{c} \(1- \frac{\beta-1}{s_0}\)} 
	\( \norm{\phi}_{\infty}^{\frac{p_{c}}{s_{0}} -1} \norm{\nabla \phi}_{2\beta}
	+ \norm{\psi}_{\infty}^{\frac{p_{c}}{s_{0}} -1} \norm{\nabla \psi}_{2\beta} \)^{ \beta -1}.
\end{align*}
Taking $p_{c}/s_{0} = \gamma \in \(1,  p_{c}/ (\beta-1)\)$, the proof is completed in $p_{c}<1$,
because the third and last terms of \eqref{nes32:1} are the same. 
The case $p_{c} \geq 1$  is easier, so we omit the proof.
\end{proof}


\section{Construction of the modified wave operators} \label{Sec:4}

In this section, we shall prove Theorem \ref{thm:1}. By the integral equation \eqref{K2/1-1}, it is easy to notice that 
\begin{align} \label{K2/2-1}
\left\| u - u_p \right\|_{X_b} \lesssim  \left\|  \int_{- \infty}^t  U(t,s) \left( F(u(s)) - F(u_p(s)) \right)  ds \right\|_{X_b}  + \left\| \CAL{K} (t) \right\|_{X_b} + \left\| \CAL{E}_r (t)  \right\|_{X_b}  ,
\end{align} 
where $ \norm{\phi}_{X_{b}}$ is defined by \eqref{K1}.
Here, in order to deduce the integrability in $(-\infty, r_0]$ for integrands of reminder terms $\CAL{K}$ and $\CAL{E}_r$, we take $\beta < \alpha$ in spite of $u_- \in H^{0, \alpha}$. This gap yields an additional $s$ decay, see Lemma \ref{lem:hl2}.  As for the first term of the RHS of \eqref{K2/2-1}, we employ the nonlinear estimates in Section {\ref{Sec:3}} and herein we give the time decay estimate for this term; however, it is, especially in $d=3$, quite difficult to obtain the bound  $ \|  \int_{- \infty}^t  U(t,s) \left( F(u(s)) - F(u_p(s)) \right)  ds \|_{X_b}  \leq (\mbox{small}) \times \| u-u_p \|_{X_b} $ since the lack of the derivativeness of nonlinearity. Indeed, as for $f,g \in H^{\alpha}$, it is difficult to show $\| e^{ i |f|^{\alpha -1}} f - e^{i |g|^{\alpha -1} } g   \|_{H^{\beta}} \lesssim \| f-g \|_{H^{\beta}}  $ for $\alpha -1 < \beta \leq \alpha$ . Because of this reason, the contraction mapping argument on $(X_{b,R}, {d}_0)$ with ${d}_0(\phi, \psi) = \| \phi - \psi \|_{X_b} $ does not work. To avoid this difficulty, we alternatively use weak norm distance $d$ and show that even if such weak norm, $(X_{b,R}, d)$ is a complete metric space.
\subsection{Estimates for the external term}
We first give estimates associated with the prescribed asymptotic profile.
\begin{Lem} 
Define $\widehat{w}$ as in \eqref{thm:1ap}. For any $|s| \geq 1$, the followings are valid:
\begin{align}
	\norm{\widehat{w}(s)}_{H^{\alpha}} \lesssim{}& \J{ \log |s|}^{\lceil \alpha \rceil} \norm{u_{-}}_{H^{0, \alpha}}\( 1 + \norm{u_{-}}_{H^{0, \alpha}}^{p_c} \)^{\lceil \alpha \rceil}, \label{pro:a1} \\
	\norm{F\(\widehat{w}(s)\)}_{H^{\alpha}} 
	\lesssim{}& \J{ \log |s|}^{\lceil \alpha \rceil} \norm{u_{-}}_{H^{0, \alpha}}^{1+p_{c}} \( 1 + \norm{u_{-}}_{H^{0, \alpha}}^{p_c} \)^{\lceil \alpha \rceil}. \label{pro:a2}
\end{align}
\end{Lem}
\begin{proof}
Noting $\norm{\widehat{u_{-}}}_{\infty} \lesssim \norm{u_{-}}_{H^{0, \alpha}}$ given by the Sobolev embedding, the estimates have been already shown in \cite[Proposition 3.1 and Proposition 3.2]{KaMi23}. 
\end{proof}

Recall that $X_{b}$- norm is defined by \eqref{K1}.
To prove Theorem \ref{thm:1}, we have to control the $X_{b}$-norm for ${A}(t)$ and $\mathcal{E}_{r}(t)$ defined as in \eqref{inte:1remi}.

\begin{Lem} \label{lem:hl1}
Take $0< b < \(1 - \beta/2 \)(1- 2 \lambda)$. Then it holds that
\begin{align*}
	\norm{ \mathcal{A}(t)}_{X_{b}}
	\lesssim{}& 
	\norm{u_{-}}_{H^{0, \beta}}^{1+p_{c}}  \( 1 + \norm{u_{-}}_{H^{0, \beta}}^{p_c} \)^{\lceil \beta \rceil}.
\end{align*}
\end{Lem}
\begin{proof}
The estimate of $\norm{\mathcal{A}(t)}_{2}$ is easier, so we only treat $\norm{U(0,t) \mathcal{A}(t)}_{H^{0, \beta}}$.
Set $f(s) = \zeta_2 (s) a_{2}^{-1} |s|^{\lambda-1}$. 
Noting that $|f(s)| \geq 1/2$ follows form \eqref{K13-2},
using \eqref{K13}, we get
\begin{align*}
	\abs{|f(s)|^{- \frac{1}{1- \lambda}} -1} \lesssim \abs{|f(s)| -1} \int_{0}^{1} \frac{d \theta}{\(1 + \theta (|f(s)| -1) \)^{\frac{1}{1-\lambda}+1}} 
	\lesssim \abs{ \frac{\zeta_2 (s)}{|s|^{1-\lambda}} - a_{2}}
	\lesssim |s|^{-(1- 2\lambda)} 
\end{align*}
for any $s \leq -r_{0}$.
Hence it follows from \eqref{pro:a2} that
\begin{align*}
	\norm{U(0,t) \mathcal{A}(t)}_{H^{0, \beta}} \leq{}& 
	\norm{\int_{-\infty}^{t}  |x|^{\beta} \mathcal{F}^{-1} \left( \frac{c_+s}{|\zeta _2 (s)|^{1/(1- \lambda)} } -1 \right) F(\widehat{w}(s)) \frac{ds}{c_+ |s|}}_{2} \\
	\lesssim{}& \int_{-\infty}^{t} \abs{|f(s)|^{- \frac{1}{1- \lambda}} -1} \norm{F(\widehat{w}(s))}_{\dot{H}^{\beta}} \frac{ds}{c_+ |s|} \\
	\lesssim{}& |t|^{-(1- 2 \lambda)} \J{ \log |t|}^{\lceil \beta \rceil} \norm{u_{-}}_{H^{0, \alpha}}^{1+p_{c}} \( 1 + \norm{u_{-}}_{H^{0, \alpha}}^{p_c} \)^{\lceil \alpha \rceil}
\end{align*}
holds for any $t \leq -r_0$. 
This completes the proof.
\end{proof}

The next statement is the estimates for $\mathcal{E}_{r}(t)$.
\begin{Lem} \label{lem:hl2}
Let $0 < b < ( \alpha- \beta)(1- 2\lambda)/2$.
Then these hold that 
\begin{align*}
	\norm{R(t) \widehat{w}(t)}_{X_{b}} \lesssim{}& \norm{u_{-}}_{H^{0, \alpha}} \( 1 + \norm{u_{-}}_{H^{0, \alpha}}^{p_c} \)^{\lceil \alpha \rceil}, \\
	\norm{\int_{-\infty}^{t} U(t,s) R(s) F(\widehat{w}(s)) \frac{ds}{c_{+} |s|}}_{X_{b}}
 \lesssim{}& \norm{u_{-}}_{H^{0, \alpha}}^{1+p_{c}}  \( 1 + \norm{u_{-}}_{H^{0, \alpha}}^{p_c} \)^{\lceil \alpha \rceil}.
\end{align*}
\end{Lem}
\begin{proof}
Applying \eqref{K13} and \eqref{le33:1},
we have
\begin{align*}
	\norm{R(t) \widehat{w}(t)}_{2} 
	={}& \norm{(\mathcal{M}_2(t) - 1) \mathcal{F}^{-1} \widehat{w}(t)}_{2}
	\lesssim |t|^{-\theta (1-2\lambda)}\norm{ \widehat{w}(t) }_{\dot{H}^{2\theta}}
\end{align*}
for any $t \leq -r_0$. Taking $\theta = \alpha/2 \in [0,1]$, one sees from \eqref{pro:a1} that
\begin{align*}
	\norm{R(t) \widehat{w}(t)}_{2} \lesssim{}& |t|^{- \frac{\alpha}{2} (1-2\lambda)}  \J{ \log |t|}^{\lceil \alpha \rceil} \norm{u_+}_{H^{0, \alpha}}\( 1 + \norm{\widehat{u_+}}_{\infty}^{p_c} \)^{\lceil \alpha \rceil}
\end{align*}
for any $t \leq -r_0$.
Arguing as in the above, by \eqref{opa:3}, it holds that
\begin{align*}
	\norm{|J(t)|^{\beta} R(t) \widehat{w}(t)}_{2} 
	\lesssim{}&  \norm{(\mathcal{M}_2(t) - 1) \mathcal{F}^{-1} |\nabla|^{\beta} \widehat{w}(t)}_{2} \\
	\lesssim{}& |t|^{- \frac{\alpha-\beta}{2} (1-2\lambda)}  \J{ \log |t|}^{\lceil \alpha \rceil} \norm{u_+}_{H^{0, \alpha}}\( 1 + \norm{\widehat{u_+}}_{\infty}^{p_c} \)^{\lceil \alpha \rceil}
\end{align*}
for any $t \leq -r_0$. Here we take $\theta = \alpha -\beta$.
Hence the first estimate is valid.
Combining \eqref{pro:a2} with \eqref{opa:3}, the second estimate similarly follows.
\end{proof}

\subsection{Estimation for leading terms}

The subsequent estimate is relevant to $F(u(s)) - F(u_p(s)$.

\begin{Lem} \label{lem:ne1}
Assume $u_{-} \in H^{0, \beta}$. 
Let $b>0$.
Then it holds that
\begin{align*}
	{}&\sup_{t \leq -r_{0}} |t|^{\frac{\beta}{2}(1-2\lambda)+b} \norm{\int_{-\infty}^{t} U(t,s) \( F(u(s)) - F(u_p(s))\)\, ds}_{2} \\
	\lesssim{}&  \( \norm{u_{-}}_{H^{0, \beta}} + \norm{u-u_p}_{X_{b}} \)^{p_c} \norm{u-u_p}_{X_{b}}.
\end{align*}
\end{Lem}
\begin{proof}
The nonlinearity can be decomposed as $F(u) - F(u_p) = F^{(1)} (u) + F^{(2)}(u)$, where
\begin{align}
	F^{(1)}(u) = \chi_{\{|u_p| \leq |u -u_p|\}} \(F(u) - F(u_p)\), \quad
	F^{(2)}(u) = \chi_{\{|u_p| \geq |u -u_p|\}} \(F(u) - F(u_p)\),
	\label{pr:31c}
\end{align}
and $\chi_A$ is a characteristic function on $A \subset \R^{1+d}$. 
Since
\[
	|F(u) - F(u_p)| \lesssim \( |u-u_p|^{p_c+1} + |u_p|^{p_c}|u-u_p| \),
\]
we have
\[
	|F^{(1)}(u)| \lesssim |u-u_p|^{p_c+1}, \quad |F^{(2)}(u)| \lesssim |u_p|^{p_c}|u-u_p|.
\]
As for $F^{(2)}(u)$, it holds that
\begin{align*}
	\norm{F^{(2)}(u)}_{2} \lesssim{}&  \norm{u_p}_{\infty}^{p_c} \norm{u-u_p}_{2} 
	\lesssim |s|^{-1-\frac{\beta}{2}(1-2\lambda)-b} \norm{\widehat{u_{-}}}_{\infty}^{p_c} \norm{u-u_p}_{X_{b}}.
\end{align*}
This implies that
\begin{align}
	\norm{\int_{-\infty}^{t} U(t,s) F^{(2)}(u(s))\, ds}_{2} 
	\lesssim{}& |t|^{-\frac{\beta}{2}(1-2\lambda)-b} \norm{\widehat{u_{-}}}_{\infty}^{p_c}  \norm{u-u_p}_{X_{b}} 
	\label{nest:1}
\end{align}
for any $t \leq -r_0$. 
Let us handle the estimate for $F^{(1)}(u)$. 
By the Gagliardo-Nirenberg inequality and \eqref{opa:1}, we obtain
\begin{align}
\begin{aligned}
	\norm{u-u_p}_{\infty} \lesssim{}&  \norm{u-u_p}_{2}^{1-\theta} \norm{\mathcal{M}_1(s) | \nabla| \mathcal{M}_1(-s)(u-u_p)}_{2}^{\theta} \\
	\lesssim{}& |s|^{-\frac1{p_c} - \frac{\beta}{2}(1-2\lambda)(1-\theta)-b} \norm{u-u_p}_{X_{b}}
\end{aligned}
	\label{kes:1}
\end{align}
for any $s \leq -r_0$, where $\theta = d/(2\beta) \in (0,1)$.
Hence since 
\begin{align*}
	\norm{F^{(1)}(u)}_{2} \lesssim{}& \norm{u-u_p}_{\infty}^{p_c} \norm{u-u_p}_{2} 
	\lesssim |s|^{-1 - \frac{\beta}{2}(1-2\lambda)-b - \frac{\beta}{2}(1-2\lambda)(1-\theta)p_c-b p_c} \norm{u-u_p}_{X_{b}}^{1+p_c},
\end{align*}
it is established that 
\begin{align*}
	\norm{\int_{-\infty}^{t} U(t,s) F^{(1)}(u(s))\, ds}_{2} 
	\lesssim{}& \norm{u-u_p}_{X_{b}}^{1+p_c} \int_{-\infty}^{t} |s|^{-1 - \frac{\beta}{2}(1-2\lambda)-b - \frac{\beta}{2}(1-2\lambda)(1-\theta)p_c-b p_c}\, ds \\
	\lesssim{}& |t|^{- \frac{\beta}{2}(1-2\lambda)-b - \frac{\beta}{2}(1-2\lambda)(1-\theta)p_c-b p_c} \norm{u-u_p}_{X_{b}}^{1+p_c} 
\end{align*}
for any $t \leq -r_0$. 
Thus there exists $b_0 = b_0(\lambda, d, \beta, \alpha)>0$ such that
\begin{align}
	|t|^{\frac{\beta}{2}(1-2\lambda)+b} \norm{\int_{-\infty}^{t} U(t,s) F^{(1)}(u(s))\, ds}_{2} 
	\lesssim{}& |t|^{- b_0} \norm{u-u_p}_{X_{b}}^{1+p_c}
	\label{nest:2}
\end{align}
for any $t \leq -r_0$.
Combining \eqref{nest:1} with \eqref{nest:2}, we conclude the desired estimate.
\end{proof}

\begin{Prop} \label{pro:31}
Assume $u_{-} \in H^{0, \beta}$. 
Then there exists $b  = b(\lambda, d, \beta) >0$ such that
\begin{align*}
	&{}\sup_{t \leq -r_{0}} |t|^{b} \norm{|J(t)|^{\beta} \int_{-\infty}^{t} U(t,s) \( F(u(s)) - F(u_p(s))\)\, ds}_{2} \\
	\lesssim{}& \(\norm{u-u_p}_{X_{b}} + \norm{\widehat{u_{-}}}_{\infty} \)^{p_c} \norm{u-u_p}_{X_{b}}
\end{align*}
is valid if $d=1$. When $d=2$, $3$, the estimate 
\begin{align*}
	&{}\sup_{t \leq -r_{0}} |t|^{b} \norm{|J(t)|^{\beta} \int_{-\infty}^{t} U(t,s) \( F(u(s)) - F(u_p(s))\)\, ds}_{2} \\
	\lesssim{}& \(\norm{u-u_p}_{X_{b}} + \norm{u_{-}}_{H^{0, \beta}} \)^{p_c} \norm{u-u_p}_{X_{b}}
	+ \norm{u_{-}}_{H^{0, \beta}} \norm{u-u_p}_{X_{b}}^{p_{c}} \\
	&{}+ \( \norm{u-u_p}_{X_{b}} + \norm{u_{-}}_{H^{0, \beta}} \)^{\gamma(\beta-1)+1} 
	\norm{u-u_p}_{X_{b}}^{p_{c} - \gamma (\beta-1)}
\end{align*}
holds for any $\gamma \in (1,  p_{c}/ (\beta-1))$ if $p_{c} <1$, otherwise
\begin{align*}
	&{}\sup_{t \leq -r_{0}} |t|^{b} \norm{|J(t)|^{\beta} \int_{-\infty}^{t} U(t,s) \( F(u(s)) - F(u_p(s))\)\, ds}_{2} \\
	\lesssim{}& \(\norm{u-u_p}_{X_{b}} + \norm{u_{-}}_{H^{0, \beta}} \)^{p_c} \norm{u-u_p}_{X_{b}} \\
	&{}+\( \norm{u-u_p}_{X_{b}} + \norm{u_{-}}_{H^{0, \beta}} \)^{p_c-1+ \beta} \norm{u-u_p}_{X_{b}}^{2- \beta}.
\end{align*}
\end{Prop}

\begin{proof}
Set $\phi = \mathcal{M}^{-1}_1(t) u$ and $\psi = \mathcal{M}_1^{-1}(t) u_p$.
we see from \eqref{opa:1} that 
\begin{align}
\begin{aligned}
	&{}\norm{|J(t)|^{\beta} \int_{-\infty}^{t} U(t,s) \( F(u(s)) - F(u_p(s))\)\, ds}_{2} \\
	\leq{}& \int_{-\infty}^{t} |\zeta_2(s)|^{\beta} \norm{|\nabla|^{\beta} \( F(\phi(s)) - F(\psi(s))\)}_{2}\, ds.
\end{aligned}
	\label{p31:1}
\end{align}
Arguing as in \eqref{kes:1}, it is obtained that 
\begin{align}
\begin{aligned}
	\norm{\phi}_{\infty} \lesssim{}& |s|^{-\frac{1}{p_c}} \(|s|^{- \beta_{\lambda}-b} \norm{u-u_p}_{X_{b}} + \norm{\widehat{u_{-}}}_{\infty} \), \\
	\norm{\phi - \psi}_{\infty} \lesssim{}& |s|^{-\frac{1}{p_c} - \beta_{\lambda} -b} \norm{u-u_p}_{X_{b}} 
\end{aligned}
	\label{p31:2}
\end{align}
for any $s \leq -r_0$, where $\beta_{\lambda} = 2^{-1} \beta (1-2\lambda)(1-\theta)$ and $\theta = d/(2\beta) \in (0,1)$. 
Also we easily verify from \eqref{pro:a1} that
\begin{align}
\begin{aligned}{}
	&{}\norm{\psi}_{\infty} \lesssim |s|^{-\frac{1}{p_c}} \norm{\widehat{u_{-}}}_{\infty}, 
	\quad \norm{\psi}_{\dot{H}^{\beta}} \lesssim |s|^{- \beta (1-\lambda)} \J{ \log |s|}^{\lceil \beta \rceil} \norm{u_{-}}_{H^{0, \beta}}
\end{aligned}
	\label{p31:3}
\end{align} 
for any $s \leq -r_0$.
We see from $\norm{\nabla \varphi}_{2\beta} \lesssim \norm{\varphi}_{\infty}^{1-\frac{1}{\beta}} \norm{\varphi}_{\dot{H}^{\beta}}^{\frac{1}{\beta}}$ and \eqref{p31:3} that 
\begin{align}
	\norm{\nabla \psi}_{2 \beta} \lesssim |s|^{- \(1 - \frac{1}{\beta} \) \frac{1}{p_{c}} - (1- \lambda) } \J{ \log |s|}^{\lceil \beta \rceil}
	\norm{\widehat{u_{-}}}_{\infty}^{1- \frac{1}{\beta}} \norm{u_{-}}_{H^{0, \beta}}^{\frac{1}{\beta}}.
	\label{825:2}
\end{align}
Let us deal with the case $d=1$. 
By means of \eqref{p31:2} and \eqref{p31:3}, we deduce from Lemma \ref{lem:34} that for any $s \leq -r_0$,
\begin{align*}
	&{}|\zeta_2(s)|^{\beta} \norm{|\nabla|^{\beta} \( F(\phi(s)) - F(\psi(s))\)}_{2} \\
	\lesssim{}&  \left(\norm{\phi - \psi}_{\infty} + \norm{\psi}_{\infty} \right)^{p_c} \norm{|\alpha(s)|^{\beta} (u-u_p)(s)}_{2} \\
	&{}+ \J{ \log |s|}^{\lceil \beta \rceil} \left(\norm{\phi - \psi}_{\infty} + \norm{\psi}_{\infty} \right)^{p_c-1} \norm{\phi-\psi}_{\infty} \norm{u_{-}}_{H^{0,\beta}} \\
	\lesssim{}& |s|^{-1-b} \left(|s|^{- \beta_{\lambda}-b} \norm{u-u_p}_{X_{b}}  + \norm{\widehat{u_{-}}}_{\infty} \right)^{p_c} \norm{u-u_p}_{X_{b}} \\
	&{}+ \J{ \log |s|}^{\lceil \beta \rceil} |s|^{-1-b- \beta_{\lambda}} \left(s^{- \beta_{\lambda}-b} \norm{u-u_p}_{X_{b}} + \norm{\widehat{u_{-}}}_{\infty} \right)^{p_c-1} 
	\norm{u_{-}}_{H^{0,\beta}} \norm{u-u_p}_{X_{b}}.
\end{align*}
Hence it follows from \eqref{p31:1} that there exists $b_1 = b_1(d, \lambda, \beta) > 0$ such that
\begin{align*}
	&{}|t|^{b} \norm{|J(t)|^{\beta} \int_{-\infty}^{t} U(t,s) \( F(u(s)) - F(u_p(s))\)\, ds}_{2} \\
	\lesssim{}& \(|t|^{- b_1-b} \norm{u-u_p}_{X_{b}} + \norm{\widehat{u_{-}}}_{\infty} \)^{p_c} \norm{u-u_p}_{X_{b}} \\
	&{}+ |t|^{- b_1} \(|t|^{- b_1-b} \norm{u-u_p}_{X_{b}} + \norm{\widehat{u_{-}}}_{\infty} \)^{p_c-1} \norm{u_{-}}_{H^{0, \beta}} \norm{u-u_p}_{X_{b}}
\end{align*}
for any $t \leq -r_0$, from which the desired assertion follows.
We shall estimate the case $d=2$, $3$.
Firstly the Gagliardo-Nirenberg inequality and \eqref{pro:a1} tell us that
\begin{align}
\begin{aligned}
	\norm{\phi}_{\dot{H}^{\beta}} \lesssim{}& |s|^{-\beta(1-\lambda)} \J{ \log |s|}^{\lceil \beta \rceil} \(|s|^{-b} \norm{u-u_p}_{X_{b}} + \norm{u_{-}}_{H^{0, \beta}} \), \\
	\norm{\phi - \psi}_{\dot{H}^{\beta}} \lesssim{}& |s|^{-\beta(1-\lambda) -b} \J{ \log |s|}^{\lceil \beta \rceil} \norm{u-u_p}_{X_{b}}. 
\end{aligned}
\label{p31:4}
\end{align}
Further it comes from \eqref{825:2} that
\begin{align}
\begin{aligned}
	\norm{\nabla \phi}_{2 \beta} \lesssim{}& |s|^{- \(1 - \frac{1}{\beta} \) \frac{1}{p_{c}} - (1- \lambda)} \J{ \log s}^{\lceil \beta \rceil} \\
	&{}\times \(|s|^{- \(1 - \frac{1}{\beta} \) \beta_{\lambda}-b} \norm{u-u_p}_{X_{b}} + \norm{\widehat{u_{-}}}_{\infty}^{1- \frac{1}{\beta}} \norm{u_{-}}_{H^{0, \beta}}^{\frac{1}{\beta}} \), \\
	\norm{\nabla \(\phi - \psi\)}_{ 2 \beta} \lesssim{}& |s|^{- \(1 - \frac{1}{\beta} \) \( \frac{1}{p_{c}} + \beta_{\lambda} \)
	 - (1- \lambda) -b } \J{ \log |s|}^{\lceil \beta \rceil} \norm{u-u_p}_{X_{b}} 
\end{aligned}
\label{825:1}
\end{align}
Arguing as in the case $d=1$, in terms of the first term, we estimate
\begin{align*}
	&{}|\zeta_2(s)|^{\beta} \norm{\phi}_{\infty}^{p_c+1-\beta} \norm{\nabla \phi}_{2\beta}^{\beta-1} \norm{\nabla (\phi - \psi)}_{2\beta} \\
	\lesssim{}& \J{ \log |s|}^{\beta \lceil \beta \rceil} |s|^{-1- b - \( 1 - \frac{1}{\beta} \) \beta_{\lambda} } 
	\( |s|^{- \beta_{\lambda}-b} \norm{u-u_p}_{X_{b}} + \norm{\widehat{u_{-}}}_{\infty} \)^{p_{c} + 1 - \beta} \\
	&{} \times \( |s|^{-\(1 - \frac{1}{\beta}\)\beta_{\lambda} - b} \norm{u-u_p}_{X_{b}} + \norm{\widehat{u_{-}}}_{\infty}^{1- \frac{1}{\beta}} \norm{u_{-}}_{H^{0, \beta}}^{\frac{1}{\beta}} \)^{\beta -1} \norm{u-u_p}_{X_{b}}
\end{align*}
for any $s \leq -r_0$.
The second term are led to
\begin{align*}
	|\zeta_2(s)|^{\beta} \norm{\phi}_{\infty}^{p_c} \norm{\phi - \psi}_{\dot{H}^{\beta}} 
	\lesssim{}& |s|^{-1 -b} \( |s|^{- \beta_{\lambda} -b} \norm{u-u_p}_{X_{b}} + \norm{\widehat{u_{-}}}_{\infty} \)^{p_c}
	\norm{u-u_p}_{X_{b}}.
\end{align*}

From now on, we shall only deal with the case $p_{c} < 1$, because the case $p_{c} \geq 1$ is similar.
In terms of the third term, one obtains 
\begin{align*}
	|\zeta_2(s)|^{\beta}\norm{\phi - \psi}_{\infty}^{p_{c}} \norm{ \psi}_{\dot{H}^{ \beta}} 
	\lesssim{}& \J{ \log |s|}^{\lceil \beta \rceil}
	|s|^{- 1- p_{c}(\beta_{\lambda} + b)} 
	\norm{u_{-}}_{H^{0, \beta}}
	\norm{u-u_p}_{X_{b}}^{p_{c}}
\end{align*}
for any $s \leq -r_0$. 
Further, as for the last term, it is calculated that
\begin{align*}
	&{}|\zeta_{2}(s)|^{ \beta} \norm{\phi - \psi}_{\infty}^{p_{c} - \gamma( \beta -1)} 
	\( \norm{\phi}_{\infty}^{\gamma -1} \norm{\nabla \phi}_{2\beta}
	+ \norm{\psi}_{\infty}^{\gamma -1} \norm{\nabla \psi}_{2\beta} \)^{ \beta -1} \norm{\nabla \psi}_{2 \beta} \\
	\lesssim{}& \J{ \log |s|}^{\beta \lceil \beta \rceil}
	|s|^{- 1- (\beta_{\lambda} +b)(p_{c} - \gamma( \beta -1))} 
	\norm{u_{-}}_{H^{0, \beta}} \( |s|^{- \beta_{\lambda}-b} \norm{u-u_p}_{X_{b}} + \norm{\widehat{u_{-}}}_{\infty} \)^{(\gamma-1)(\beta-1)} \\
	&{} \quad \times \( |s|^{-\(1 - \frac{1}{\beta}\)\beta_{\lambda} - b} \norm{u-u_p}_{X_{b}} + \norm{\widehat{u_{-}}}_{\infty}^{1- \frac{1}{\beta}} \norm{u_{-}}_{H^{0, \beta}}^{\frac{1}{\beta}} \)^{\beta -1}
	\norm{u-u_p}_{X_{b}}^{p_{c} - \gamma( \beta -1)} 
\end{align*}
for any $s \leq -r_0$, as long as
$0 < b < \beta_{\lambda} \frac{p_{c}- \gamma ( \beta-1)}{1- p_{c}+ \gamma ( \beta-1)}$.
Connecting these above, we conclude that there exists $b_1 = b_1(d, \lambda, \beta) > 0$ such that
\begin{align*}
	&{}|t|^{b} \norm{|J(t)|^{\beta} \int_{-\infty}^{t} U(t,s) \( F(u(s)) - F(u_p(s))\)\, ds}_{2} \\
	\lesssim{}& |t|^{-b_1} \(|t|^{-b_{1}-b} \norm{u-u_p}_{X_{b}} + \norm{\widehat{u_{-}}}_{\infty} + \norm{u_{-}}_{H^{0, \beta}} \)^{p_c} \norm{u-u_p}_{X_{b}} \\
	&{}+ \( |t|^{-b_{1}-b} \norm{u-u_p}_{X_{b}} + \norm{\widehat{u_{-}}}_{\infty} \)^{p_c}
	\norm{u-u_p}_{X_{b}} \\
	&{}+ |t|^{- b_1} 	\norm{u_{-}}_{H^{0, \beta}} \norm{u-u_p}_{X_{b}}^{p_{c}} \\
	&{}+ |t|^{-b_1} \( |t|^{-b_{1}-b} \norm{u-u_p}_{X_{b}} + \norm{\widehat{u_{-}}}_{\infty} \)^{(\gamma-1)(\beta-1)} \\
	&{}\quad \times \( |t|^{-b_{1}-b} \norm{u-u_p}_{X_{b}} + \norm{\widehat{u_{-}}}_{\infty} + \norm{u_{-}}_{H^{0, \beta}} \)^{\beta-1} 
	\norm{u_{-}}_{H^{0, \beta}} \norm{u-u_p}_{X_{b}}^{p_{c} - \gamma (\beta-1)}
\end{align*}
for any $t \leq -r_0$. 
By the Sobolev embedding, this yields the desired estimate.
\end{proof}

\subsection{Proof of Theorem \ref{thm:1}}
We are now in position to prove Theorem \ref{thm:1}.

Take $b>0$ given by Proposition \ref{pro:31}.
From the integral equation \eqref{inte:1}, we define the map 
\begin{align*} 
\Psi(v(t)) = i \int_{-\infty}^{t} U(t,s) \left( F(v(t) + u_{p}(s)) - F(u_p(s)) \right)  ds + \CAL{A}(t) + \CAL{E}_{\mathrm{r}} (t)
\end{align*}
for any $v \in X_{b, R}$.
We remark that $(X_{b, R}, d)$ is the complete metric space.
Indeed, any Cauchy sequence $\{v_n\} \subset X_{b, R}$ has a strong limit $v \in L^{\infty}(- \infty, -r_{0}; L^2)$ and 
the boundedness of $X_{b, R}$ provides us a subsequence $\{|J(t)|^{\beta} v_{n_k}(t)\}$ and $w(t) \in L^2$ such that $|J(t)|^{\beta}v_{n_k}(t) \rightarrow w(t)$ in $L^2$ for each $t \geq T$. 
Since $|J(t)|^{ \beta}$ is the closed operator in $\D{|J (t)|^{\beta}} \coloneqq \{ u \in L^2 \mid \| |J (t)|^{\beta} u \|_{2} < \infty \}$, in view of the uniqueness of the limit in $\mathscr{S}'$, it holds that $w(t) = |J(t)|^{\beta}v(t)$ in $\mathscr{S}'$ and hence
\begin{align*}
	&{}t^{\frac{\beta}{2}(1-2\lambda) +b} \norm{v(t)}_{2} + t^b \norm{|J(t)|^{\beta} v(t)}_{2} \\
	\leq{}& \liminf_{k \rightarrow \infty} \{t^{\frac{\beta}{2}(1-2\lambda) +b} \norm{v_{n_k}(t)}_{2} 
	+ t^b  \norm{|J(t)|^{\beta} v_{n_k}(t)}_{2} \} \leq R 
\end{align*}
for any $t \geq T$. Thus $v \in X_{b, R}$. 
As for the external term, 
by Lemma \ref{lem:hl1} and Lemma \ref{lem:hl2}, 
\begin{align*}
	\norm{\CAL{A}( \cdot) + \CAL{E}_{\mathrm{r}} (\cdot)}_{X_{b}}
	\leq C \(1 + \norm{u_{-}}_{H^{0, \alpha}} \)^{3p_{c}} \norm{u_{-}}_{H^{0, \alpha}}
	\leq C_{0} \norm{u_{-}}_{H^{0, \alpha}}
\end{align*}
if $ \varepsilon_{1} \leq 1$.
Let us here show that $\Psi \colon X_{b, R} \rightarrow X_{b, R}$ 
is a contraction map with $R = 2C_{0} \norm{u_{-}}_{H^{0, \alpha}}$. 
We only deal with the case $p_c<1$ because the other case is similar.
Recall $d/2 < \beta < \alpha < \min(d, 1+p_{c})$.
Combining Lemma \ref{lem:ne1} with Proposition \ref{pro:31}, we see that
there exists $C>0$ such that
\begin{align*}
	\norm{\Psi(v)}_{X_{b}} \leq{}& C\(  R+ \norm{u_{-}}_{H^{0, \alpha}} \)^{p_c} R
	+C \norm{u_{-}}_{H^{0, \alpha}} R^{p_{c}} 
	+C \(R+ \norm{u_{-}}_{H^{0, \alpha}} \)^{\gamma(\beta-1)+1} R^{p_{c} - \gamma (\beta-1)} \\
	&{}+C_{0} \norm{u_{-}}_{H^{0, \alpha}} \\
	\leq{}& b_{1} \norm{u_{-}}_{H^{0, \alpha}}^{p_{c}+1} +  C_{0} \norm{u_{-}}_{H^{0, \alpha}}
\end{align*}
for any $v \in X_{T, b, R}$, as long as $\varepsilon_{1} \leq 1$.
Hence, we have $\norm{\Psi(v)}_{X_{b}} \leq R$ if $b_{1} \epsilon_{1}^{p_{c}} \leq C_{0}$.
Thus $\Psi (v) \in X_{b,R}$ holds. 
In the same manner,  it holds that 
\begin{align*}
	d\(\Psi(v_1), \Psi(v_2)\) 
	\leq{}& a_{1} \norm{u_{-}}_{H^{0, \alpha}}^{p_{c}} \norm{v_{1} - v_{2}}_{L^{ \infty} (- \infty, -r_{0} ; L^{2}) } 
	\leq \frac12 d(v_1, v_2)
\end{align*}
for any $v_1$, $v_2 \in X_{b, R}$ whenever $a_{1} \varepsilon_{1}^{p_{c}} \leq 1/2$.
These imply that $\Psi \colon X_{b, R} \to X_{b, R}$ is the contraction map. 
Thus, we have a unique solution in $X_{b, R}$.
This completes the proof.

\subsection{Proof of Corollary \ref{cor:waveop} }
We finish this section by proving the existence of the wave operator. Fix $\varepsilon \in (0, \varepsilon_1]$, where $ \varepsilon_{1}$ is given in Theorem \ref{thm:1}.
We see from Theorem \ref{thm:1} that for any $u_{-} \in H^{0, \alpha}$ with $\norm{u_{-}}_{H^{0, \alpha}} \leq \varepsilon$, 
\eqref{eq1} admits a unique $\mathcal{F}H^{\beta}$-solution $u$ on $(- \infty, -r_0]$
satisfying
\[
	\norm{U(0, t)\left(u(t)- u_{p}(t)\right)}_{H^{0, \beta}} 
	\leq C \varepsilon |t|^{-\mu}
\]
for all $t \leq - r_0$.
By means of \eqref{pro:a1} and the Sobolev embedding,
this implies 
\begin{align*}
	\norm{U(0, -r_{0}) u(-r_{0})}_{H^{0, \beta}} 
	\leq{}& \norm{U(0, -r_{0})u_{p}(-r_{0})}_{H^{0, \beta}} + C \varepsilon |r_{0}|^{-\mu} \\
	\leq{}& \norm{u_{-}}_{2} + C(\eta, a_{2}, r_{0}) \norm{u_{-}}_{H^{0, \beta}} \( 1 + \norm{u_{-}}_{H^{0, \beta}}^{p_c} \)^{\lceil \alpha \rceil}
	+ C(r_{0}, \mu) \varepsilon \\
	\leq{}& b_{1} \varepsilon
\end{align*}
and thus $U(0, -r_{0}) u(-r_{0}) \in H^{0, \beta}$.
By Proposition \ref{pro:ext1}, we have the solution on $(- \infty, 0]$ satisfying
\begin{align*}
	\norm{u(0)}_{H^{0, \beta}}
	\leq{}& \norm{|J( \cdot)|^{ \beta}u}_{L^{ \infty}(-r_{0}, 0; L^{2})} 
	\leq A_{0} \norm{U(0, -r_{0}) u(-r_{0})}_{H^{0, \beta}}
	\leq A_{0}b_{1} \varepsilon
\end{align*}
where $A_{0}>0$ depends on 
$r_{0}$ and $ \varepsilon_{1}$. 
Hence, setting $c_{0} = A_{0}b_{1}$, we can define a map by 
\begin{align*}
	\mathcal{W}_{-} \colon B^{\alpha}_{\varepsilon} \ni u_{-} \mapsto u(0) \in B^{ \beta}_{c_{0} \varepsilon}
\end{align*}
for any $ \varepsilon \in (0, \varepsilon_{1}]$. This is the so-called modified wave operator.

\section{Construction of the modified inverse wave operator} \label{Sec:5}

Thanks to Corollary \ref{cor:waveop} and Proposition \ref{pro:ext1}, we have the propagator from $u_-$ to $u(r_0)$, and in this section we construct the modified inverse operator which propagates from $u(r_0)$ to $u_{+}$. 
Let us first introduce the function space
\begin{align*} 
	&\mathscr{X}_{}^{ \varepsilon, A} \coloneqq \{ u \in C([r_0, \infty) ; L^{2} \cap L^{ \infty} (\R^{d})) \mid 
	\norm{u}_{\mathscr{X}_{}^{ \varepsilon, A}} < \infty \} , \\ 
	& \norm{u}_{\mathscr{X}_{}^{\varepsilon, A} } 
	\coloneqq{} \sup_{t \geq r_0} |t|^{-A \varepsilon^{\frac{1}{d(1-\lambda)}}} 
	\norm{\J{J(t)}^{ \beta} u}_{2}
	+ \sup_{t \geq r_0} |t|^{\frac{d}{2}(1- \lambda)} \norm{u(t)}_{\infty}
\end{align*}
for some small $ \varepsilon>0$ and a certain $A>0$. 
The smallness of $ \varepsilon$ are needed to construct the global solution on $[r_{0}, \infty)$.
As for the detail, see the beginning of Appendix \ref{app:A}.
As for the constant $A$, we refer to comments above Proposition \ref{prop4.5} in Appendix \ref{app:A}. 
Here we set 
\begin{align}
\begin{aligned}
v(t) &{}= \CAL{M}_+ U(0,t) u(t), \\ 
\SCR{A} (t) &{} =
\begin{cases}
0, & d=1, \\ 
  |\widehat{v}(t,x)|^{p_c} - \left( |t|^{- \chi} + |\widehat{v}(t,x)  |^2  \right)^{\frac{p_c}2}   , & d=2,3, 
 \end{cases}
  \\ 
	I_1(t) &{}= \mathcal{F} \( \mathcal{M}_{2}^{-1} - 1 \) \mathcal{F}^{-1} 
	\( |\mathcal{F} \mathcal{M}_{2} v|^{p_{c}} \mathcal{F} \mathcal{M}_{2} v \),
	\\[2mm]
	I_2(t) &{}=  |\mathcal{F} \mathcal{M}_{2} v|^{p_{c}} \mathcal{F} \mathcal{M}_{2} v - |\mathcal{F} v|^{p_{c}} \mathcal{F} v,
\end{aligned}
	\label{remain:1}
\end{align} 
and
\[
	\theta(t) = \theta( t, x) =
	\begin{cases} \vspace{0.1cm}
	  {\displaystyle\eta \int_{r_{0}}^{t}  |\widehat{v}(s, x)|^{p_c} |\zeta_{2}(s)|^{- \frac{1}{1- \lambda}}\, ds }, & d=1, \\ 
	 {\displaystyle \eta \int_{r_{0}}^{t} \( |s|^{- \chi} + |\widehat{v}(s, x)|^{2} \)^{\frac{p_{c}}{2}} |\zeta_{2}(s)|^{- \frac{1}{1- \lambda}}\, ds, }& d=2,3.
	 \end{cases} 
\]
We further define 
\begin{align*}
	 \mathcal{P}(t) {}:=  \exp (i \theta(t)), \quad v_{p} \coloneqq{}& 
	 \mathcal{F}^{-1} \mathcal{P}(t) \mathcal{F}v =\mathcal{F}^{-1} e^{i \theta (t)}  \mathcal{F}v.
\end{align*}
Then by the similar arguments in Section {4} of \cite{HN06}, we have 
\begin{align}
\begin{aligned}
	 \widehat{v_{p}}(t) - \widehat{v_{p}}(s) ={}& 
	- i \eta \int_{s}^{t} |\zeta_{2}(\tau)|^{-\frac{1}{1- \lambda}} 
	\mathcal{P}(\tau) \mathscr{A}(\tau)
	\widehat{v}(\tau)\, d \tau \\
	&{}- i \eta \int_{s}^{t} \mathcal{P}(\tau) \( I_1(\tau)+I_2(\tau) \) |\zeta_{2}(\tau)|^{- \frac{1}{1- \lambda}}\, d \tau.
\end{aligned}
	\label{iv:12}
\end{align}

In this section, we focus on the existence of inverse modified wave operators $\CAL{W}_+ ^{-1}: u(0) \mapsto u_+ $, namely, we establish that $\widehat{v_{p}}(t) \in H^{\delta}$ has a limit by proving $ \widehat{v_{p}}(t) - \widehat{v_{p}}(s) \rightarrow 0$ as $t$, $s \rightarrow \infty$ on $H^{\delta}$ under global solutions $u \in \mathscr{X}^{\ep, A}$ to \eqref{eq1}.
We first need to show that \eqref{eq1} with $u_0 = u(r_0)$ is global well-posed. However, as mentioned in Section \ref{ssec:1}, we introduce the result on GWP as Theorem \ref{iv:thm1}  and give a sketch of the proof in Appendix \ref{app:A}, since result on GWP has been essentially proven in \cite{KM21}.
Hence we only treat that $\widehat{v_{p}}(t)$ has a limit in $H^{\delta}$ in this section.

In the following, we handle the RHS of \eqref{iv:12} in $H^{\delta}$. If one can obtain the integrand of which are integrable in $[r_0, \infty)$, we immediately find that ${v_p}(t)$ has the limit in $H^{0,\delta}$, that is exactly the $u_+$. Hence the main purpose of this section is to seek the decay estimate of the integrand of the RHS of \eqref{iv:12}, namely we show the following Propositions:

\begin{Prop} \label{pro:gam1}
Let $d = 2$, $3$.
Assume $d/2 < \delta < \beta < \min(1+p_c, d)$.
Let $u$ be the solution given by Theorem \ref{iv:thm1}.
Set $ \epsilon_{ \lambda} = A \varepsilon^{\frac{1}{d(1-\lambda)}}$.
Then it holds that
\begin{align*}
	\norm{\mathcal{P}(\tau) \mathscr{A}(\tau) \widehat{v}(\tau)}_{H^{\delta}} 
	\lesssim{}&  \varepsilon \tau^{-\frac{\chi}{2} \kappa + 2 \epsilon_{\lambda}},
\end{align*}
where $\kappa = \kappa (d, \lambda, \delta)$ satisfies $\kappa = 1 + p_{c} - \delta$ if $p_{c} \geq 1$ and 
$\kappa =\min \( 3 p_{c} - \delta, (1- p_{c})(1+ p_{c} - \delta) + p_{c}^2 \)$
if $p_{c} <1$.
\end{Prop}

\begin{Prop} \label{Pro:56}
Let $d \leq 3$.
Assume $d/2 < \delta < \beta < \min(1+p_c, d)$.
Let $u$ be the solution given by Theorem \ref{iv:thm1}.
Set $ \epsilon_{ \lambda} = A \varepsilon^{\frac{1}{d(1-\lambda)}}$.
If $d=1$, then
\begin{align*}
	\norm{\mathcal{P}(\tau)\(I_{1}(\tau) + I_{2}(\tau) \) }_{H^{\delta}} \lesssim \varepsilon \tau^{-\frac{ \beta - \delta}{2}(1- 2 \lambda) +  \(p_{c} + 1\) \epsilon_{ \lambda} }
\end{align*}
holds. When $d=2$, $3$, 
the following is valid:
\begin{align*}
	\norm{ \mathcal{P}(\tau)\(I_{1}(\tau) + I_{2}(\tau) \)}_{H^{\delta}}
	\lesssim{}& \varepsilon \tau^{- \frac{ \beta - \delta}{2}(1 + \min (1, p_{c}) - \delta)(1 - 2 \lambda) + \frac{\chi}{2} +  (2p_{c} + 1) \epsilon_{ \lambda}},
\end{align*}
where $I_{j}(\tau)$ are defined by \eqref{remain:1}.
\end{Prop}

With these propositions in place, let us initially give the proof of Theorem \ref{thm:iv2}.
We postpone the proof of the propositions to the end of this section. 

\subsection{Proof of Theorem \ref{thm:iv2}}
Let us prove Theorem \ref{thm:iv2} by Proposition \ref{pro:gam1} and Proposition \ref{Pro:56}.
We shall only treat the case $p_{c} < 1$, because the others are similar.
Set 
\[
	0 < \chi < (\beta - \delta)(1 + p_{c} - \delta)(1 - 2 \lambda).
\]
We see from Proposition \ref{pro:gam1} and Proposition \ref{Pro:56} that 
there exists $\mu_{0} > 0$ satisfying 
\begin{align*}
	0 < \mu_{0} < \min \( \frac{\chi}{2} \kappa, \frac{ \beta - \delta}{2}(1 + p_{c} - \delta)(1 - 2 \lambda) - \frac{\chi}{2} \),
\end{align*}
such that
\begin{align*}
	\norm{ \widehat{v_{p}}(t) -\widehat{ v_{p}}(s)}_{H^{\delta}}
	\lesssim{}& 
	\int_{s}^{t} 
	\norm{\mathcal{P}(\tau) \mathscr{A}(\tau)
	\widehat{v}(\tau)}_{H^{\delta}}|\zeta_{2}(\tau)|^{-\frac{1}{1- \lambda}} \, d \tau \\
	&{}+ \int_{s}^{t} \norm{\mathcal{P}(\tau) \( I_1(\tau)+I_2(\tau) \)}_{H^{\delta}} |\zeta_{2}(\tau)|^{- \frac{1}{1- \lambda}}\, d \tau \\
	\lesssim{}& \varepsilon \( t^{- \mu_{0} + (2p_{c} + 1) \epsilon_{ \lambda}} - s^{-\mu_{0} + (2p_{c} + 1) \epsilon_{ \lambda}} \)
\end{align*}
for any $r_{0} < s<t$.
Then there exists a unique limit $u_{+} \in H^{0, \delta}$ such that
\begin{align*}
	\norm{{v_{p} }(t) - u_{+}}_{H^{0,\delta}} \lesssim{}& \varepsilon t^{- \mu_{0} + (2p_{c} + 1) \epsilon_{ \lambda}}
\end{align*}
for any $t > r_{0}$. 
This completes the proof of Theorem \ref{thm:iv2}.


\subsection{Proof of Proposition \ref{pro:gam1} and Proposition \ref{Pro:56}.}
Before showing these propositions, we introduce three important lemmata. 

By the Gagliardo-Nirenberg inequality and \eqref{le33:1}, 
we easily verify the followings for $\mathcal{M}_{+} U(0, t)u$ frequently used afterward.
\begin{Lem} \label{leps:1}
Let $u$ be the solution given by Theorem \ref{iv:thm1} and  $v = \mathcal{M}_{+} U(0, t)u$. 
Set $ \epsilon_{ \lambda} = A \varepsilon^{\frac{1}{d(1-\lambda)}}$.
Fix $ \alpha \in [0, \beta]$.
Then the following hold: For any $t \geq r_{0}$,
\begin{itemize}
\item $\norm{\widehat{v}(t)}_{\dot{H}^{\alpha}} + \norm{\mathcal{F} \mathcal{M}_{2} v}_{\dot{H}^{\alpha}} 
	\lesssim \varepsilon |t|^{ \frac{ \alpha}{ \beta} \epsilon_{ \lambda}}$, \;
	$\norm{\mathcal{F} \mathcal{M}_{2} v}_{\infty} + \norm{ \widehat{v}}_{\infty}
	\lesssim \varepsilon |t|^{ \frac{d}{ 2\beta} \epsilon_{ \lambda}}$, \;
	$\norm{\nabla \widehat{v}(t)}_{2 \delta} \lesssim \varepsilon t^{\epsilon_{\lambda}}$.
\item $\norm{\mathcal{F} \(\mathcal{M}_{2} - 1\) v}_{\infty} 
	\lesssim \varepsilon |t|^{- \alpha ( 1- 2 \lambda) + \epsilon_{ \lambda}}$, 
	\;
	$\norm{ \mathcal{F} \(\mathcal{M}_{2} - 1\) v}_{\dot{H}^{ \alpha}} 
	\lesssim \varepsilon |t|^{- \frac{ \beta - \alpha}{2} ( 1- 2 \lambda) + \epsilon_{ \lambda}}$,
\end{itemize}
where $\alpha \in [0, ( \beta - d/2)/2)$. 
\end{Lem}

In the modified wave operator, we need a modification of the phase correction $e^{i \theta (t)}$ to remove a non-integrable part for $t$.
Moreover, in $d=2$, $3$, one employs a further modification such as $\mathcal{P}(t)$ which is firstly used in \cite{GO93}, due to the lack of differentiability of the nonlinearity.
Hence the following sharp inequalities is needed to extract necessary time decay for a term arising from the modification.
Set $\chi >0$, $n \in \N \cup \{0\}$ and $t \geq r_{0}$.
We denote 
\begin{align}
	\mathscr{A}_{n}(t) ={}& \mathscr{A}_{n}(t ; z) = \( |z|^{p_{c}-n} - \(|t|^{- \chi} + |z|^{2} \)^{\frac{p_{c}-n}{2}} \) z^{n},
	\label{scra:1}
\end{align}
for any $z \in \C$. Further, define 
\begin{align*}
	\mathscr{B}_{n}(t ; z) = \(|t|^{- \chi} + |z|^{2} \)^{\frac{p_{c}-n}{2}} z^{n}, \quad
	\widetilde{\mathscr{B}_{n}}(t ; z) = \(|t|^{- \chi} + |z|^{2} \)^{\frac{p_{c}-n}{2}} |z|^{n} 
\end{align*}
for any $z \in \C$, where we remark that 
$$ \mathscr{A}(t) = \mathscr{A}_{0}(t ; \hat{v}(t)), \quad d=2,3. $$

\begin{Lem} \label{lem:vv:1}
The following estimates hold: 
\begin{itemize}
\item If $p_{c} \in (0,2)$, then 
\begin{align}
	&{} \abs{\mathscr{A}_{n}(t ; z)} \lesssim |t|^{- \frac{p_{c}}{2}\chi }. \label{mph:1}
\end{align}
\item When $p_{c} \in (0,1)$, 
\begin{align}
	&{} \abs{\widetilde{\mathscr{B}_{n}}(t ; z) - \widetilde{\mathscr{B}_{n}}(t ; w)} 
	+ \abs{\mathscr{B}_{n}(t ; z) - \mathscr{B}_{n}(t ; w)} \lesssim |z_{ }- w{}|^{p_{c}}, \label{mph:2} \\
	&{} \abs{\mathscr{A}_0(t ; z) z} \leq |t|^{- p_{c} \chi} |z|^{1- p_{c}}. \label{mph:3} 
\end{align}
\end{itemize}
 
\end{Lem}

\begin{proof}
We shall prove \eqref{mph:1} by an induction argument, with reference to the proof of \cite[Lemma 2.4]{AIMMU}.
The case $n=0$ immediately follows from the concavity of $|z|^{p_{c}/2}$.
When $n \in \N$, one computes
\begin{align*}
	\mathscr{A}_{n+1}(t ; z)
	={}& \left\{ \( |z|^{p_{c}-n} - (|\tau|^{- \chi} + |z|^{2})^{\frac{p_{c}-n}{2}} \) z^{n} \right\} |z|^{-1} z \\
	&{}+ (|\tau|^{- \chi} + |z|^{2})^{\frac{p_{c}-(n+1)}{2}} z^{n+1} |z|^{-1} \left\{ (|\tau|^{- \chi} + |z|^{2})^{\frac{1}{2}} -  |z| \right\}.
\end{align*}
Thanks to the hypothesis of the induction, the first term is bounded by $C |t|^{- \frac{p_{c}}{2}\chi}$.
The second term is also dominated by $C |t|^{- \frac{p_{c}}{2}\chi}$ from the direct calculation.
This completes the proof of \eqref{mph:1}.
We shall next prove \eqref{mph:2} by the induction.
We only deal with $\abs{\mathscr{B}_{n}(t ; z) - \mathscr{B}_{n}(t ; w)} \lesssim |z - w|^{p_{c}}$, because the other part is similar.
By the concavity of $|z|^{p_{c}}$, one has
\begin{align*}
	\abs{ \mathscr{B}_{0}(t; z) - \mathscr{B}_{0}(t ; w)}
	\leq{}& \abs{ \frac{|z|^{2} - |w|^{2}}{(|\tau|^{- \chi} + |z|^{2})^{\frac{1}{2}} + (|\tau|^{- \chi} + |w|^{2})^{\frac{1}{2}}}}^{p_{c}} 
	\leq \abs{z-w}^{p_{c}},
\end{align*}
which is the case $n=0$.
When $n \in \N$, we compute
\begin{align*}
	&{}\mathscr{B}_{n+1}(t; z) - \mathscr{B}_{n+1}(t ; w) \\
	={}& \left\{ \(|t|^{- \chi} + |z|^{2} \)^{\frac{p_{c}-n}{2}} z^{n} - \(|t|^{- \chi} + |w|^{2} \)^{\frac{p_{c}-n}{2}} w^{n} \right\}  
	\(|t|^{- \chi} + |z|^{2} \)^{-\frac{1}{2}} z \\
	&{}+ \(|t|^{- \chi} + |w|^{2} \)^{\frac{p_{c}-n}{2}} w^{n} 
	\left\{ \(|t|^{- \chi} + |z|^{2} \)^{-\frac{1}{2}} z 
	- \(|t|^{- \chi} + |w|^{2} \)^{-\frac{1}{2}} w \right\}.
\end{align*}
The first term is bounded by $C|z-w|^{p_{c}}$ from the hypothesis of the induction.  
If $|z-w| \geq \(|t|^{- \chi} + |w|^{2} \)^{\frac{1}{2}}$, then the second term is also dominated by $C |z-w|^{p_{c}}$.
When $|z-w| \leq \(|t|^{- \chi} + |w|^{2} \)^{\frac{1}{2}}$, by a further calculation,
\begin{align*}
	&{}\(|t|^{- \chi} + |w|^{2} \)^{\frac{p_{c}-n}{2}} w^{n}
	\left\{ \(|t|^{- \chi} + |z|^{2} \)^{-\frac{1}{2}} z - \(|t|^{- \chi} + |w|^{2} \)^{-\frac{1}{2}} w \right\} \\
	={}& \(|t|^{- \chi} + |w|^{2} \)^{\frac{p_{c}-(n+1)}{2}} w^{n} \(|t|^{- \chi} + |z|^{2} \)^{-\frac{1}{2}} 
	\left\{ \(|t|^{- \chi} + |w|^{2} \)^{\frac{1}{2}} - \(|t|^{- \chi} + |z|^{2} \)^{\frac{1}{2}} \right\} z \\
	&{}+ \(|t|^{- \chi} + |w|^{2} \)^{\frac{p_{c}-(n+1)}{2}} w^{n} (z-w).
\end{align*}
Hence, the second term is bounded by $C|z -w|^{p_{c}}$. The proof  of \eqref{mph:2} is completed.
Let us finally show \eqref{mph:3}.
By means of the concavity of $|z|^{p_{c}}$, we obtain
\begin{align*}
	\abs{\mathscr{A}_0(t ; z) z}
	\leq{}& \abs{ |z| - \(|t|^{- \chi} + |z|^{2} \)^{\frac{1}{2}} }^{p_{c}} |z| 
	\leq \abs{\frac{|t|^{- \chi}}{|z| + \(|t|^{- \chi} + |z|^{2} \)^{\frac{1}{2}}} }^{p_{c}} |z|,
\end{align*}
as desired.
\end{proof}

\begin{Lem} \label{lem:mp1}
Let $d \leq 3$.
Assume $d/2 < \delta < \beta < \min(1+p_c, d)$.
Let $u$ be the solution given by Theorem \ref{iv:thm1}.
Set $ \epsilon_{ \lambda} = A \varepsilon^{\frac{1}{d(1-\lambda)}}$.
Then when $d=1$, it holds that
\begin{align*}
	\norm{\mathcal{P}(\tau)}_{\dot{H}^{\delta}} 
	\lesssim{}& \tau^{p_{c} \epsilon_{\lambda}}.
\end{align*}
Under $d=2$, $3$, 
if $p_{c} < 1$, then
\begin{align*}
	\norm{\mathcal{P}(\tau)}_{\dot{H}^{\delta}} 
	\lesssim{}& \varepsilon \tau^{\frac{\chi}{2}(1- p_{c}) + \frac{\chi}{2}(\delta -1) + \delta \epsilon_{\lambda}}
\end{align*}
holds, otherwise, 
\begin{align*}
	\norm{\mathcal{P}(\tau)}_{\dot{H}^{\delta}}  
	\lesssim \varepsilon \tau^{\frac{\chi}{2}(2- p_{c})(\delta-1) + p_{c} \epsilon_{\lambda}}.
\end{align*}
\end{Lem}

\begin{proof}[\textit{\bf Proof of Lemma \ref{lem:mp1}}]
When $d=1$, noting $| \mathcal{P}( \tau, x) -  \mathcal{P}(\tau, y)| \lesssim |\theta(\tau, x) - \theta(\tau, y)|$ and $p_{c} \geq 2$, 
we see from a use of the representation via the Gagliardo semi-norm that
\begin{align*}
	\norm{\mathcal{P}( \tau)}_{\dot{H}^{ \delta}}
	\lesssim{}& \norm{ \theta(\tau)}_{\dot{H}^{ \delta}} 
	\lesssim \int_{r_{0}}^{\tau} \norm{\( |s|^{- \chi} + |\widehat{v}(s)|^{2} \)^{\frac{p_{c}-2}{2}}}_{L^{\infty}} 
	\norm{\widehat{v}(s)}_{\dot{H}^{ \delta}}
	|\zeta_{2}(s)|^{- \frac{1}{1- \lambda}}\, ds 
	\lesssim \varepsilon \tau^{p_{c} \epsilon_{ \lambda}}.
\end{align*}
Let us turn to the case $d=2$, $3$. We only handle the case $p_{c} < 1$, because the other case is easier.
In view of $\nabla \mathcal{P}(\tau) = i \nabla \theta(\tau) \mathcal{P}(\tau)$, 
we see from Lemma \ref{lem_lei} that
\begin{align*}
	\norm{\mathcal{P}( \tau)}_{\dot{H}^{ \delta}}
	= \norm{ \nabla \mathcal{P}( \tau)}_{\dot{H}^{\delta-1}}
	\lesssim{}& \norm{\theta(\tau)}_{\dot{H}^{\delta}} + \norm{\nabla \theta(\tau)}_{2 \delta} \norm{|\nabla|^{\delta-1} \mathcal{P}(\tau)}_{\frac{2\delta}{\delta-1}}  
\end{align*}
The Gagliardo-Nirenberg inequality gives us 
\begin{align}
	\norm{|\nabla|^{\delta-1} \mathcal{P}(\tau)}_{\frac{2\delta}{\delta-1}} 
	\lesssim{}& \norm{ \mathcal{P}(\tau)}_{\infty}^{2- \delta} \norm{\nabla \mathcal{P}(\tau)}_{2\delta}^{\delta-1} 
	= \norm{\nabla \theta(\tau)}_{2 \delta}^{\delta-1},
	\label{96:6}
\end{align}
which yields
\begin{align}
	\norm{\mathcal{P}( \tau)}_{\dot{H}^{ \delta}}
	\lesssim{}& \norm{\theta(\tau)}_{\dot{H}^{\delta}} + \norm{\nabla \theta(\tau)}_{2 \delta}^{ \delta}.
	\label{96:1}
\end{align}
Since
\begin{align*}
	\nabla \theta(\tau) = \eta \int_{r_{0}}^{\tau} p_{c} \( |s|^{- \chi} + |\widehat{v}(s)|^{2} \)^{\frac{p_{c}-2}{2}} 
	\Re \( \overline{\widehat{v}(s)} \nabla \widehat{v}(s) \) |\zeta_{2}(s)|^{- \frac{1}{1- \lambda}}\, ds,
\end{align*}
we see from Lemma \ref{leps:1} that
\begin{align}
\begin{aligned}
	\norm{\nabla \theta(\tau)}_{2 \delta} 
	\lesssim{}& \int_{r_{0}}^{\tau} \norm{\( |s|^{- \chi} + |\widehat{v}(s)|^{2} \)^{\frac{p_{c}-1}{2}}}_{\infty} 
	\norm{\nabla \widehat{v}(s)}_{2 \delta} |\zeta_{2}(s)|^{- \frac{1}{1- \lambda}}\, ds \\
	\lesssim{}& \int_{r_{0}}^{\tau} |s|^{-1 + \frac{\chi}{2}(1 - p_{c})} \norm{\nabla \widehat{v}(s)}_{2 \delta}\, ds
	\lesssim \varepsilon \tau^{\frac{\chi}{2}(1- p_{c}) + \epsilon_{\lambda}}. 
\end{aligned}
	\label{96:2}
\end{align}
Further, applying Lemma \ref{lem_lei}, one has
\begin{align}
\begin{aligned}
	&{}\norm{\( |s|^{- \chi} + |\widehat{v}(s)|^{2} \)^{\frac{p_{c}-2}{2}} 
	\overline{\widehat{v}(s)} \nabla \widehat{v}(s)}_{\dot{H}^{\delta-1}} \\
	\lesssim{}& \norm{|\nabla|^{\delta-1} \( |s|^{- \chi} + |\widehat{v}(s)|^{2} \)^{\frac{p_{c}-2}{2}} 
	\overline{\widehat{v}(s)}}_{\frac{2\delta}{\delta-1}} \norm{ \nabla \widehat{v}(s)}_{2 \delta} \\
	&{}+ \norm{\( |s|^{- \chi} + |\widehat{v}(s)|^{2} \)^{\frac{p_{c}-2}{2}} 
	\overline{\widehat{v}(s)}}_{\infty} \norm{\widehat{v}(s)}_{\dot{H}^{\delta}}.
\end{aligned}
	\label{96:3}
\end{align}
By means of the Gagliardo-Nirenberg inequality, we see that 
\begin{align*}
	&{}\norm{|\nabla|^{\delta-1} \( |s|^{- \chi} + |\widehat{v}(s)|^{2} \)^{\frac{p_{c}-2}{2}} 
	\overline{\widehat{v}(s)}}_{\frac{2\delta}{\delta-1}} \\
	\lesssim{}& \norm{\( |s|^{- \chi} + |\widehat{v}(s)|^{2} \)^{\frac{p_{c}-2}{2}} 
	\overline{\widehat{v}(s)}}_{\infty}^{2- \delta} 
	 \norm{ |\nabla| \( |s|^{- \chi} + |\widehat{v}(s)|^{2} \)^{\frac{p_{c}-2}{2}} 
	\overline{\widehat{v}(s)} }_{2\delta}^{\delta-1} \\
	\lesssim{}& s^{\frac{\chi}{2}(1- p_{c}) + \frac{\chi}{2}(\delta -1)} 
	\norm{ \nabla \widehat{v}(s)}_{2\delta}^{\delta-1}.
\end{align*}
Hence by Lemma \ref{leps:1}, one obtains 
\begin{align*}
	\norm{\( |s|^{- \chi} + |\widehat{v}(s)|^{2} \)^{\frac{p_{c}-2}{2}} 
	\overline{\widehat{v}(s)} \nabla \widehat{v}(s)}_{\dot{H}^{\delta-1}} 
	\lesssim{}& s^{\frac{\chi}{2}(1- p_{c}) + \frac{\chi}{2}(\delta -1)} \norm{ \nabla \widehat{v}(s)}_{2\delta}^{\delta}
	+\varepsilon s^{\frac{\chi}{2}(1- p_{c}) + \epsilon_{\lambda}} \\
	\lesssim{}& \varepsilon s^{\frac{\chi}{2}(1- p_{c}) + \frac{\chi}{2}(\delta -1) + \delta \epsilon_{\lambda}}.
\end{align*}
This yields
\begin{align}
\begin{aligned}
	\norm{\theta(\tau)}_{\dot{H}^{\delta}} 
	\lesssim{}& \int_{r_{0}}^{\tau} \norm{\( |s|^{- \chi} + |\widehat{v}(s)|^{2} \)^{\frac{p_{c}-2}{2}} 
	\overline{\widehat{v}(s)} \nabla \widehat{v}(s)}_{\dot{H}^{\delta-1}} |\zeta_{2}(s)|^{- \frac{1}{1- \lambda}}\, ds \\
	\lesssim{}& \varepsilon \tau^{\frac{\chi}{2}(1- p_{c}) + \frac{\chi}{2}(\delta -1) + \delta \epsilon_{\lambda}}.
	\label{96:5}
\end{aligned}
\end{align}
Plugging \eqref{96:5} and \eqref{96:2} into \eqref{96:1}, we conclude
\begin{align*}
	\norm{\mathcal{P}( \tau)}_{\dot{H}^{ \delta}}
	\lesssim{}& \varepsilon \tau^{\frac{\chi}{2}(1- p_{c}) + \frac{\chi}{2}(\delta -1) + \delta \epsilon_{\lambda}}
	+ \varepsilon \tau^{\frac{\chi}{2}(1- p_{c}) \delta + \delta \epsilon_{\lambda}} 
	\lesssim \varepsilon \tau^{\frac{\chi}{2}(1- p_{c}) + \frac{\chi}{2}(\delta -1) + \delta \epsilon_{\lambda}},
\end{align*}
as desired. 
This completes the proof.
\end{proof}

Using these lemmata, we show Proposition \ref{pro:gam1} and Proposition \ref{Pro:56}.

\begin{proof}[\textit{\bf Proof of Proposition \ref{pro:gam1}}]
From \eqref{mph:1} and Lemma \ref{leps:1}, we easily show
\begin{align*}
	\norm{\mathcal{P}(\tau) \mathscr{A}(\tau) \widehat{v}(\tau)}_{2} 
	\leq \norm{ \mathscr{A}( \tau)}_{\infty} \norm{\widehat{v}(\tau)}_{2} 
	\lesssim \varepsilon \tau^{- \frac{p_{c}}{2} + \epsilon_{ \lambda}}.
\end{align*}
Also, the direct computation shows
\begin{align*}
	\nabla \( \mathcal{P}(\tau) \mathscr{A}(\tau) \widehat{v}(\tau) \)
	={}& \nabla \mathcal{P}(\tau) \mathscr{A}(\tau) \widehat{v}(\tau)
	+ \frac{p_{c}}{2} \mathcal{P}(\tau) |\mathscr{A}_{2}(\tau)| \nabla \widehat{v}(\tau) \\
	&{}+ \frac{p_{c}}{2} \mathcal{P}(\tau) \mathscr{A}_{2}(\tau) \nabla \widehat{v}(\tau)
	+ \mathcal{P}(\tau) \mathscr{A}(\tau) \nabla \widehat{v}(\tau),
\end{align*}
where recall that $ \mathscr{A}_{2}( \tau)$ is defined by \eqref{scra:1}. 
This implies 
\begin{align}
\begin{aligned}
	\norm{\mathcal{P}(\tau) \mathscr{A}(\tau) \widehat{v}(\tau)}_{\dot{H}^{\delta}}
	\lesssim{}& \norm{\nabla \mathcal{P}(\tau) \mathscr{A}(\tau) \widehat{v}(\tau)}_{\dot{H}^{\delta-1}}
	+ \norm{\mathcal{P}(\tau) |\mathscr{A}_{2}(\tau)| \nabla \widehat{v}(\tau)}_{\dot{H}^{\delta-1}} \\
	&{}+ \norm{\mathcal{P}(\tau) \mathscr{A}_{2}(\tau) \nabla \widehat{v}(\tau)}_{\dot{H}^{\delta-1}}
	+ \norm{\mathcal{P}(\tau) \mathscr{A}(\tau) \nabla \widehat{v}(\tau)}_{\dot{H}^{\delta-1}} \\
	=: &{} J_{1}(\tau) + J_{2}(\tau) + J_{3}(\tau) + J_{4}(\tau).
\end{aligned}
	\label{827:6}
\end{align}
In view of $\nabla \mathcal{P}(\tau) = i \nabla \theta(\tau) \mathcal{P}(\tau)$, 
we see from Lemma \ref{lem_lei} that
\begin{align}
\begin{aligned}
	J_{1} (\tau)
	\lesssim{}& \norm{ \mathcal{P}(\tau)}_{\dot{H}^{\delta}}  \norm{\mathscr{A}(\tau) \widehat{v}(\tau)}_{\infty} 
	+ \norm{\nabla \theta(\tau)}_{2 \delta} \norm{|\nabla|^{\delta-1} \mathscr{A}(\tau) \widehat{v}(\tau)}_{\frac{2\delta}{\delta-1}}. 
\end{aligned}
	\label{826:2}
\end{align}
Hereafter, one shall address the case $p_{c} <1$.
Collecting \eqref{mph:1} and \eqref{mph:3}, the Gagliardo-Nirenberg inequality gives us
\begin{align*}
	\norm{|\nabla|^{\delta-1} \mathscr{A}(\tau)\widehat{v}(\tau)}_{\frac{2\delta}{\delta-1}} 
	\lesssim{}&\norm{\mathscr{A}(\tau)\widehat{v}(\tau)}_{\infty}^{2- \delta} \norm{\nabla \( \mathscr{A}(\tau)\widehat{v}(\tau) \)}_{2\delta}^{\delta-1} \\
	\lesssim{}& \norm{\mathscr{A}(\tau)\widehat{v}(\tau)}_{\infty}^{2- \delta} 
	\( \norm{\mathscr{A}_{2}(\tau)}_{\infty} + \norm{\mathscr{A}(\tau)}_{\infty} \)^{\delta -1} \norm{\nabla \widehat{v}(\tau)}_{2 \delta}^{\delta-1} \\
	\lesssim{}& \tau^{-\frac{\chi}{2} p_{c} (3 -\delta)} \norm{\widehat{v}}_{\infty}^{(1-p_{c})(2-\delta)} \norm{\nabla \widehat{v}(\tau)}_{2 \delta}^{\delta-1}.
\end{align*}
Noting \eqref{mph:3}, substituting these above into \eqref{826:2}, applying Lemma \ref{lem:mp1}, we establish
\begin{align}
\begin{aligned}
	J_{1} (\tau) \lesssim{}& \varepsilon \tau^{\frac{\chi}{2}(1- p_{c}) + \frac{\chi}{2}(\delta -1) + \delta \epsilon_{\lambda}} (\varepsilon \tau^{\epsilon_{\lambda}})^{(1-p_{c})} \tau^{-p_{c} \chi} 
	+ \varepsilon \tau^{\frac{1 - p_{c}}{2} \chi - \frac{3-\delta}{2} p_{c} \chi  + 2\epsilon_{\lambda}} \\
	\lesssim{}& \varepsilon \tau^{- \frac{\chi}{2}(3 p_{c} - \delta) + 2 \epsilon_{\lambda}}. 
\end{aligned}
	\label{827:1}
\end{align}
In terms of $J_{2}$, arguing as in $J_{1}$, together with Lemma \ref{leps:1}, it follows from \eqref{96:2} and \eqref{96:6} that
\begin{align*}
	J_{2} (\tau)  \lesssim{}&\norm{|\nabla|^{\delta-1} \mathcal{P}(\tau)}_{\frac{2\delta}{\delta-1}} \norm{\mathscr{A}_{2}(\tau)}_{\infty}
	\norm{\nabla \widehat{v}(\tau)}_{2 \delta} \\
	&{}+ 
	\norm{|\nabla|^{\delta-1} |\mathscr{A}_{2}(\tau)|}_{\frac{2\delta}{\delta-1}}
	\norm{\nabla \widehat{v}(\tau)}_{2 \delta} 
	+
	 \norm{\mathscr{A}_{2}(\tau)}_{\infty}
	\norm{\widehat{v}(\tau)}_{\dot{H}^{\delta}} \\
	\lesssim{}& \varepsilon \tau^{\frac{\chi}{2}(1- p_{c})(\delta-1) -\frac{\chi}{2}p_{c} + \delta \epsilon_{\lambda}} 
	+\norm{|\nabla|^{\delta-1} |\mathscr{A}_{2}(\tau)|}_{\frac{2\delta}{\delta-1}} \varepsilon \tau^{\epsilon_{\lambda}} 
	+ \varepsilon \tau^{-\frac{\chi}{2} p_{c} + \epsilon_{\lambda}}.
\end{align*}
Note that by \eqref{nl:2} and \eqref{mph:2}, $|\mathscr{A}_{2}(\tau)|$ is a $p_{c}$-H\"older continous function, namely
\[
	\abs{ |\mathscr{A}_{2}(\tau; z_1)| - |\mathscr{A}_{2}(\tau; z_{2})|}
	\lesssim  |z_1 - z_{2}|^{p_{c}}.
\]
The Gagliardo-Nirenberg inequality leads to 
\begin{align}
\begin{aligned}
	\norm{|\nabla|^{\delta-1} |\mathscr{A}_{2}(\tau)|}_{\frac{2\delta}{\delta-1}} 
	\lesssim{}& \norm{\mathscr{A}_{2}(\tau)}_{\infty}^{1- \frac{\delta-1}{s_0}} 
	\norm{ |\nabla|^{s_0}  |\mathscr{A}_{2}(\tau)| }_{\frac{2\delta}{s_0}}^{\frac{\delta-1}{s_0}},
\end{aligned}
	\label{827:a}
\end{align}
where $s_0 = p_{c} - \varepsilon_{1}$ for any $\varepsilon_{1} \in (0, p_{c} - \delta+1)$.
Applying Lemma \ref{lem_hcha}, we deduce that 
\begin{align*}
	\norm{|\nabla|^{s_0} |\mathscr{A}_{2}(\tau)|}_{\frac{2\delta}{s_0}} 
	\lesssim{}& \norm{\widehat{v}(\tau)}_{\infty}^{p_{c} - \frac{s_0}{s}} \norm{|\nabla|^s \widehat{v}(\tau)}_{\frac{2\delta}{s}}^{\frac{s_0}{s}}
\end{align*}
for any $s \in \(s_0/p_{c}, 1\)$.
One then sees from the Gagliaro-Nirenberg inequality that 
\[
	\norm{|\nabla|^s \widehat{v}(\tau)}_{\frac{2\delta}{s}} 
	\lesssim \norm{\widehat{v}(\tau)}_{\infty}^{1-s} \norm{\nabla \widehat{v}(\tau)}_{2\delta}^s,
\]
which implies
\begin{align*}
	\norm{|\nabla|^{s_0} |\mathscr{A}_{2}(\tau)|}_{\frac{2\delta}{s_0}} 
	\lesssim{}&\norm{\widehat{v}(\tau)}_{\infty}^{p_{c} -s_0} \norm{\nabla \widehat{v}(\tau)}_{2\delta}^{s_0}.
\end{align*}
Plugging the above into \eqref{827:a}, we have
\begin{align*}
	\norm{|\nabla|^{\delta-1} |\mathscr{A}_{2}(\tau)|}_{\frac{2\delta}{\delta-1}} 
	\lesssim{}& \tau^{- \frac{\chi}{2}p_{c} \(1- \frac{\delta-1}{s_0}\)} 
	\norm{\widehat{v}(\tau)}_{\infty}^{\(\frac{p_{c}}{s_{0}} -1\)(\delta-1)} \norm{\nabla \widehat{v}(\tau)}_{2\delta}^{\delta -1} \\
	\lesssim{}& \varepsilon^{\frac{p_{c}}{s_{0}}(\delta-1)} \tau^{- \frac{\chi}{2}\(p_{c} - \frac{p_{c}}{s_{0}} (\delta-1) \) + \epsilon_{\lambda}},
\end{align*}
which implies
\begin{align}
\begin{aligned}
	J_{2} (\tau) \lesssim{}& \varepsilon \tau^{\frac{\chi}{2}(1- p_{c})(\delta-1) -\frac{\chi}{2}p_{c} + \delta \epsilon_{\lambda}} 
	+ \varepsilon^{1 + \frac{p_{c}}{s_{0}}(\delta-1)} \tau^{- \frac{\chi}{2}\(p_{c} - \frac{p_{c}}{s_{0}} (\delta-1) \) + 2\epsilon_{\lambda}} 
	+ \varepsilon \tau^{-\frac{\chi}{2} p_{c} + \epsilon_{\lambda}} \\
	\lesssim{}& \varepsilon \tau^{-\frac{\chi}{2} (1- p_{c})(1+ p_{c} - \delta) -\frac{\chi}{2}p_{c}^2 + 2 \epsilon_{\lambda}}.
\end{aligned}
	\label{827:4}
\end{align}
Similarly to $J_{2}$, it holds that
\begin{align}
	J_{3} (\tau)+ J_{4} (\tau) \lesssim{}& \varepsilon \tau^{-\frac{\chi}{2} (1- p_{c})(1+ p_{c} - \delta) -\frac{\chi}{2}p_{c}^2 + 2 \epsilon_{\lambda}}.
	\label{827:5}
\end{align}
Plugging \eqref{827:1}, \eqref{827:4} and \eqref{827:5} into \eqref{827:6}, we conclude
\begin{align*}
	\norm{\mathcal{P}(\tau) \mathscr{A}(\tau) \widehat{v}(\tau)}_{\dot{H}^{\delta}} 
	\lesssim{}& \varepsilon \tau^{- \frac{\chi}{2}(3 p_{c} - \delta) + 2 \epsilon_{\lambda}} + \varepsilon \tau^{-\frac{\chi}{2} (1- p_{c})(1+ p_{c} - \delta) -\frac{\chi}{2}p_{c}^2 + 2 \epsilon_{\lambda}} 
	\lesssim \varepsilon \tau^{- \frac{\chi}{2} \mu + 2 \epsilon_{\lambda}}
\end{align*}
as long as $\kappa = \kappa (d, \lambda, \delta)$ satisfies $0 <\kappa < \min \( 3 p_{c} - \delta, (1- p_{c})(1+ p_{c} - \delta) + p_{c}^2 \)$, as desired in $p_{c}<1$. 
The case $p_{c} \geq 1$ is similar, so we omit the proof. 
\end{proof}

\begin{proof}[\textit{\bf Proof of Proposition \ref{Pro:56}}]
As in the preceding proposition, we will only treat the case $p_{c} < 1$ in $d=3$.
Let us firstly handle $I_{1}(\tau)$. 
Together with 
$\abs{\mathcal{M}^{-1}_{2}(t) - 1} \lesssim |x|^{2\theta}|t|^{- \theta (1- 2\lambda)}$ for any $ \theta \in [0,1]$, it comes from Lemma \ref{lem:34} and Lemma \ref{leps:1} that 
\begin{align}
\begin{aligned}
	\norm{I_1(\tau)}_{\dot{H}^{ \delta}}
	\leq{}& \norm{ |x|^{\delta} \( \mathcal{M}_{2}^{-1} - 1 \) \mathcal{F}^{-1} 
	\( |\mathcal{F} \mathcal{M}_{2} v|^{p_{c}} \mathcal{F} \mathcal{M}_{2} v \)}_{2} \\
	\lesssim{}& \tau^{-\frac{ \beta - \delta}{2}(1 - 2 \lambda) } 
	\norm{\mathcal{F} \mathcal{M}_{2} v}_{\infty}^{p_{c}} 
	\norm{\mathcal{F} \mathcal{M}_{2} v}_{\dot{H}^{\beta}} \\
	\lesssim{}& \varepsilon  
	\tau^{-\frac{ \beta - \delta}{2}(1- 2 \lambda) + \( p_{c} + 1 \) \epsilon_{ \lambda} }.
\end{aligned}
	\label{96:a1}
\end{align}
Also, the estimate
\begin{align}
	\norm{ \mathcal{P}(\tau) I_1(\tau) }_{2} 
	\lesssim{}& \varepsilon  
	\tau^{-\frac{\beta}{2}(1- 2 \lambda) + \( p_{c} + 1 \) \epsilon_{\lambda} }
	\label{pro:56:0}
\end{align}
is valid in the same manner.
Further, arguing as in \eqref{96:a1}, we see from the Sobolev embedding that
\begin{align}
\begin{aligned}
	\norm{I_1(\tau)}_{\infty} 
	\lesssim{}& \tau^{-\alpha(1- 2 \lambda)} \norm{\mathcal{F} \mathcal{M}_{2} v}_{\infty}^{p_{c}} 
	\norm{\mathcal{F} \mathcal{M}_{2} v}_{\dot{H}^{\beta}} \\
	\lesssim{}& \varepsilon \tau^{- \frac{1}{2} \( \beta - \frac{d}{2} \) (1 -2 \lambda) + \(  p_{c} + 1 \) \epsilon_{ \lambda} }.
\end{aligned}
	\label{96:a2}
\end{align}
for any $\alpha \in [0, ( \beta - d/2)/2)$.
By Lemma \ref{lem_lei}, we obtain
\begin{align*}
	\norm{ \mathcal{P}(\tau) I_1(\tau)}_{\dot{H}^{\delta}}
	\lesssim{}& \norm{I_1(\tau)}_{\dot{H}^{\delta}}
	+ \norm{|\nabla|^{\delta-1} \mathcal{P}(\tau)}_{\frac{2 \delta}{\delta-1}}\norm{\nabla I_1(\tau)}_{2 \delta}  \\
	&{}+\norm{\nabla \mathcal{P}(\tau)}_{2 \delta} \norm{|\nabla|^{ \delta-1} I_1(\tau)}_{\frac{2 \delta}{\delta-1}} 
	+ \norm{\mathcal{P}(\tau)}_{\dot{H}^{ \delta}} \norm{I_1(\tau)}_{\infty} .
\end{align*}
By the Gagliardo-Nireberg inequality, it follows from \eqref{96:a1} and \eqref{96:a2} that
\begin{align*}
	\norm{ \nabla I_{1}( \tau)}_{2 \delta} 
	\lesssim{}& \norm{I_{1}(\tau)}_{\infty}^{1 - \frac{1}{\delta}} \norm{I_{1}(\tau)}_{\dot{H}^{\delta}}^{\frac{1}{\delta}} 
	\lesssim \varepsilon \tau^{-\frac{ \beta - \delta}{2}(1- 2 \lambda) + \( p_{c} + 1 \) \epsilon_{ \lambda} },
\end{align*}
and
\begin{align*}
	\norm{|\nabla|^{ \delta-1} I_1(\tau)}_{L^{\frac{2 \delta}{\delta-1}}}
	\lesssim{}& \norm{I_{1}(\tau)}_{\infty}^{\frac{1}{\delta}} \norm{I_{1}(\tau)}_{\dot{H}^{\delta}}^{1- \frac{1}{\delta}} 
	\lesssim \varepsilon \tau^{-\frac{ \beta - \delta}{2}(1- 2 \lambda) + \( p_{c} + 1 \) \epsilon_{ \lambda} }.
\end{align*}
Also, in view of \eqref{96:6} and \eqref{96:2}, 
we obtain
\begin{align*}
	\norm{|\nabla|^{ \delta-1} \mathcal{P}(\tau)}_{L^{\frac{2 \delta}{\delta-1}}}
	\lesssim \varepsilon \tau^{\frac{\chi}{2}(1- p_{c})(\delta-1) + (\delta-1) \epsilon_{\lambda}}.
\end{align*}
Hence, choosing $\chi \in (0, ( \beta - \delta)( 1 - 2 \lambda))$, we establish
\begin{align*}
	\norm{ \mathcal{P}(\tau)I_1(\tau)}_{\dot{H}^{\delta}}
	\lesssim{}&  \varepsilon \tau^{- \frac{ \beta - \delta}{2}(1- 2 \lambda) + \frac{\chi}{2} + 4 \epsilon_{ \lambda} }.
\end{align*}

Let us next consider $I_{2}(\tau)$. 
We see from Lemma \ref{leps:1} and \eqref{nl:1} that
\begin{align}
\begin{aligned}
	\norm{\mathcal{P}(\tau)I_{2}(\tau)  }_{2} \lesssim{}&  
	\norm{|\mathcal{F} \mathcal{M}_{2} v|^{p_{c}} \mathcal{F} \mathcal{M}_{2} v - |\mathcal{F} v|^{p_{c}} \mathcal{F} v}_{2} \\
	\lesssim{}&  
	\( \norm{\mathcal{F} \mathcal{M}_{2} v}_{\infty} + \norm{\mathcal{F} v}_{\infty} \)^{p_{c}} \norm{\mathcal{F} \(\mathcal{M}_{2} - 1\) v}_{2} \\
	\lesssim{}& \varepsilon \tau^{- \frac{ \beta}{2} ( 1- 2 \lambda) + (p_{c} + 1) \epsilon_{ \lambda} }.
\end{aligned}
	\label{pro:56:3}
\end{align}
Also, the similar calculation gives us 
\begin{align}
	\begin{aligned}
	\norm{I_2( \tau)}_{\infty} 
	\lesssim{}& \( \norm{\mathcal{F}\mathcal{M}_{2}v}_{\infty} + \norm{\mathcal{F}v}_{\infty} \)^{p_{c}} 
	\norm{\mathcal{F} \( \mathcal{M}_{2} -1 \) v}_{\infty} \\
	\lesssim{}& \varepsilon \tau^{- \frac{1}{2} \( \beta - \frac{d}{2} \)(1- 2\lambda) + (p_{c} + 1) \epsilon_{ \lambda}}.
	\end{aligned}
	\label{iv:11}
\end{align}
We deduce from Lemma \ref{lem_lei} that
\begin{align*}
	\norm{ \mathcal{P}(\tau)I_{2}(\tau)}_{\dot{H}^{\delta}}
	\lesssim{}& \norm{I_{2}(\tau)}_{\dot{H}^{\delta}}
	+  \norm{|\nabla|^{\delta-1} \mathcal{P}(\tau)}_{\frac{2 \delta}{\delta-1}}\norm{\nabla I_2(\tau)}_{2 \delta} \\
	&{}+ \norm{\nabla \mathcal{P}(\tau)}_{2 \delta}\norm{|\nabla|^{ \delta-1} I_{2}(\tau)}_{\frac{2 \delta}{\delta-1}} 
	+ \norm{\mathcal{P}(\tau)}_{\dot{H}^{ \delta}} \norm{I_{2}(\tau)}_{\infty} .
\end{align*}

Hereafter we only handle the case $p_{c} <1$.
By Lemma \ref{lem:34} and Lemma \ref{leps:1}, one obtains 
\begin{align*}
	\norm{I_{2}( \tau)}_{\dot{H}^\delta} 
	\lesssim{}& \varepsilon \tau^{- \frac{ \beta - \delta}{2}(1 + p_{c} - \delta)(1 - 2 \lambda) +  (p_{c} + 1) \epsilon_{ \lambda}},
\end{align*}
where $\gamma \in (1,  p_{c}/ (\delta-1))$.
By the Gagliardo-Nireberg inequality, we have
\begin{align*}
	\norm{|\nabla|^{ \delta-1} I_{2}(\tau)}_{L^{\frac{2 \delta}{\delta-1}}} + \norm{ \nabla I_{2}( \tau)}_{2 \delta} 
	\lesssim \varepsilon \tau^{- \frac{ \beta - \delta}{2}(1 + p_{c} - \delta)(1 - 2 \lambda) +  (p_{c} + 1) \epsilon_{ \lambda}}.
\end{align*}
Combining these above, we conclude that
\begin{align*}
	\norm{  \mathcal{P}(\tau) I_{2}(\tau)}_{\dot{H}^{\delta}}
	\lesssim{}& \varepsilon \tau^{- \frac{ \beta - \delta}{2}( 1 + p_{c} - \delta)(1 - 2 \lambda) + \frac{\chi}{2} +  (2p_{c} + 1) \epsilon_{ \lambda}}.
\end{align*}
This completes the proof.
\end{proof}


\subsection{Proof of Corollary \ref{cor:inwa}}
We herein construct the inverse modified wave operator. Let $u_{0} \in B_{ \varepsilon}^{ \beta}$ for any $\varepsilon \in (0, \widetilde{\varepsilon}_{2}]$, where $\widetilde{\varepsilon}_{2}$ is chosen later.
Proposition \ref{pro:ext1} tells us that \eqref{eq1} admits a unique solution $u \in C([0, r_{0}]; L^{2})$ with $u(0) = u_{0}$ satisfying
\begin{align*}
	\norm{ U(0, t) u(t)}_{H^{0, \beta}} 
	\leq{}& A_{0} \norm{u_{0}}_{H^{0, \beta}}
	\leq A_{0} \varepsilon
\end{align*}
for any $t \in [0, r_{0}]$, where $A_{0} >0$ depends on $r_{0}$ and $\widetilde{\varepsilon}_{2}$.
This yields 
\begin{align*}
	\norm{ U(0, r_{0}) u(r_{0})}_{H^{0, \beta}} 	\leq A_{0} \varepsilon.
\end{align*}
Hence, taking $\widetilde{\varepsilon}_{2}>0$ such that $A_{0} \widetilde{\varepsilon}_{2} \leq \varepsilon_{2}$, 
we can apply Theorem \ref{thm:iv2} to construct a unique solution $u \in C([r_{0}, \infty); L^{2})$ with $U(0, \cdot) u \in C([r_{0}, \infty ; H^{0, \beta})$, where $\varepsilon_{1}$ is given in Theorem \ref{thm:iv2}. 
Moreover, there exists a unique $u_{+} \in H^{0, \delta}$ such that \eqref{thm:iv2a} holds.
This implies that a map
\begin{align*}
	\mathcal{W}_{+}^{-1} \colon B^{ \beta}_{ \varepsilon} \ni u(0)  \mapsto u_{+} \in H^{0, \gamma}
\end{align*}
can be defined for any $ \varepsilon \in (0, \widetilde{\varepsilon}_{2}]$.
The map is the so-called inverse wave operator.

\subsection{Construction of the modified scattering operator}

We are finally ready to define the modified scattering operator.
\begin{proof}[\textit{\bf Proof of Corollary \ref{cor:sca}}]

Take $u_{-} \in B^{ \alpha}_{ \varepsilon}$ for any $ \varepsilon \in (0, \varepsilon_{0}]$, where $\varepsilon_{0} \in (0, \epsilon_{1})$ is chosen later.
By Corollary \ref{cor:waveop}, we have the modified wave operator
\begin{align*}
	\mathcal{W}_{-} \colon B^{\alpha}_{\varepsilon} \ni u_{-}  \mapsto u(0) \in B^{ \beta}_{c_{0} \varepsilon},
\end{align*}
and thus $u(0) = W_{-} u_{-} \in B^{ \beta}_{c_{0} \varepsilon}$.
Choosing $\varepsilon_{0} >0$ as $c_{0} \varepsilon_{0} < \widetilde{\varepsilon}_{2}$, 
one sees from Corollary \ref{cor:inwa} that there exists $u_{+} \in H^{0,  \gamma}$ such that
the modified inverse wave operator $ \mathcal{W}_{+}^{-1}$ maps $u(0)$ into $u_{+}$. 
We then conclude $u_{+} = \mathcal{W}_{+}^{-1} W_{-} u_{-}$.
Therefore, the modified scattering operator $S_{+} \coloneqq \mathcal{W}_{+}^{-1} W_{-}$ can be defined on $B^{ \alpha}_{ \varepsilon}$.
This completes the proof.
\end{proof}


\section*{Appendix}

\appendix
\section{Global-in-time well-posedness on $[r_0, \infty)$} \label{app:A}
In this appendix, we shall show the following GWP result on $[r_0, \infty)$. 
\begin{Thm}[Global existence for initial value problem] \label{iv:thm1}
Let $d/2 < \beta < \min(1+p_c, d)$.
There exists $\varepsilon_2 > 0$ such that the following assertion holds: For any $\varepsilon \in (0, \varepsilon_2]$ and $u_0 \in L^{2}$ satisfying $U(0, r_{0}) u_{0} \in H^{0, \beta}$ with $\norm{U(0, r_{0})u_{0}}_{H^{0, \beta}} \leq \varepsilon$, there exists a unique $\mathcal{F}H^{\beta}$-solution to \eqref{eq1} on $[r_{0}, \infty)$ under $u(r_{0})=u_0$ satisfying $\norm{u(t)}_{\infty} \lesssim \varepsilon |t|^{-\frac{d}{2}(1- \lambda)}$
for any $t \geq r_{0}$.
\end{Thm}
To show this theorem we introduce the function space
\begin{align*}
	\mathscr{X}_{T}^{ \varepsilon, A} = \{ u \in C(I_{T}; L^{2} \cap L^{ \infty} (\R^{d})) \mid 
	\norm{u}_{\mathscr{X}_{T}^{ \varepsilon, A}} < \infty \}
\end{align*}
equipped with the distance function $d(u, v) = \norm{u-v}_{L^{ \infty}(I_{T}, L^{2})}$,
where $\norm{\, \cdot \,}_{\mathscr{X}}$ is defined by 
\begin{align*}
	\norm{u}_{\mathscr{X}} ={}& \norm{u}_{\mathscr{X}_{T}^{\varepsilon, A} } 
	\coloneqq{} \sup_{t \in I_{T}} |t|^{-A \varepsilon^{\frac{1}{d(1-\lambda)}}} 
	\norm{\J{J(t)}^{ \beta} u}_{2}
	+ \sup_{t \in I_{T}} |t|^{\frac{d}{2}(1- \lambda)} \norm{u(t)}_{\infty}
\end{align*}
for any $\varepsilon>0$ and a certain $A>0$, where $I_{T} = [r_{0}, r_{0}+T]$. 
Afterward, the smallness of $ \varepsilon$ will be imposed to extend the solution globally in-time (see Corollary \ref{cor4.4} and Proposition \ref{prop4.5}).
As for the constant $A$, we refer to comments above Proposition \ref{prop4.5}. 
Arguing as in the proof of Theorem \ref{thm:1}, it can be proven that the metric space $( \mathscr{X}_{T}^{ \varepsilon, A}, d)$ is complete. 

We will begin with the local existence of solutions to \eqref{eq1} in $\mathscr{X}_{T}^{\varepsilon,A}$.

\begin{Prop} \label{iv:1}
Let $ \varepsilon >0$, $A>0$, and $K>0$. Assume that $u_0 \in L^{2}(\R^{d})$ satisfies $| \alpha(r_{0})|^{ \beta} u_{0} \in L^{2}$
with 
$\norm{u_0}_{2} + \norm{|J(t_{0})|^{ \beta} u_{0}}_{2} \leq K$
. Then there exists $T = T(K) >0$ not depending on $ \varepsilon$ and $A$ such that \eqref{eq1} has a unique solution $u \in \mathscr{X}_{T}^{\varepsilon,A}$ with $u(r_{0}) = u_0$. Moreover the solution satisfies 
\begin{align}
	\norm{u}_{\mathscr{X}_{T}^{\varepsilon,A}} \leq B_0 K \label{iv:2}
\end{align}
for some constants $B_0 >0$ not depending on $r_0$, $T$, $A$ and $\varepsilon$.
\end{Prop}

\begin{proof}
The proof is based on the standard contraction mapping argument from the completeness of the metric space $( \mathscr{X}_{T}^{ \varepsilon, A}, d)$.
So we omit the proof.
\end{proof}
With Proposition \ref{iv:1} in place, we have the following:
\begin{Cor}
\label{cor4.4}
Let $A>0$. Then, there exists $\varepsilon_1>0$ such that the following assertion holds: 
For any $\varepsilon \in (0, \varepsilon_1]$ and $u_0 \in L^{2}$ satisfying $| \alpha(r_{0})|^{ \beta} u_0 \in L^{2}$ with $\norm{u_0}_{2} + \norm{| \alpha(r_{0})|^{ \beta} u_0}_{2} \leq \varepsilon$, there exists $T = T(\varepsilon) >0$ not depending on $A$ such that \eqref{eq1} has a unique solution $u \in \mathscr{X}_{T}^{\varepsilon, A}$ with $u(r_0) = u_0$. 
Moreover the solution satisfies 
\begin{align*}
	\norm{u}_{\mathscr{X}_{T}^{\varepsilon, A}} < \varepsilon^{\frac{1}{2}}.
\end{align*}
\end{Cor}

\begin{proof}
Taking $\varepsilon_1>0$ such that $B_0 < \varepsilon_1^{-1/2}$, we have $B_0\varepsilon < \varepsilon^{1/2}$ for any $\varepsilon \in (0,\varepsilon_1]$, where $B_0$ is the constant as in Proposition \ref{iv:1}. Applying Proposition \ref{iv:1} as $K=\varepsilon$, we obtain the desired assertion. 
\end{proof}

In what follows, we fix $A =c_0$, where $c_0$ is given in Proposition \ref{iv:pro1}. 
We show the following estimate.

\begin{Prop}[Bootstrap estimate]
\label{prop4.5}
There exists $\varepsilon_1 >0$ such that the following holds: Let $T>0$.  If $\varepsilon \in (0, \varepsilon_1]$, $u_0 \in L^{2}$ satisfies $|J(r_{0})|^{ \beta} u_0 \in L^{2}$ with $\norm{u_0}_{2} + \norm{| J(r_{0})|^{ \beta} u_0}_{2} \leq \varepsilon$, and the solution $u \in \mathscr{X}_{T}^{\varepsilon,A}$ of \eqref{eq1}  with $u(r_0) = u_0$ satisfies 
\begin{align}
\label{eq4.15}
	\norm{u}_{\mathscr{X}_{T}^{\varepsilon,A}} \leq \varepsilon^{\frac{1}{2}},
\end{align}
then there exists a constant $b_{1}$ independent of $T$ and $\varepsilon$ such that
\begin{align*}
	\norm{u}_{\mathscr{X}_{T}^{\varepsilon,A}} \leq b_{1} \varepsilon.
\end{align*}
\end{Prop}
\begin{proof}[\textit{\bf Proof of Theorem \ref{iv:thm1}}]
Once we obtain Proposition \ref{prop4.5}, Theorem \ref{iv:thm1} can be proven by the standard continuation argument (e.g. \cite[Theorem2.2]{AIMMU22}).
\end{proof}

To prove Proposition \ref{prop4.5}, we need the following lemma immediately proven by \eqref{mdfm}:

\begin{Lem} \label{iv:3}
Let $u \in C([r_{0}, \infty); L^{2})$ satisfy $|J( \cdot)|^{\beta}u \in C([r_{0}, \infty); L^{2})$ and $\alpha \in [0, ( \beta - d/2)/2)$. Then, it holds that
\begin{align*}
	\norm{u(t)}_{\infty} \lesssim |t|^{-\frac{d}{2}(1 -\lambda)} \norm{\mathcal{F} \mathcal{M}_{+} U(0, t) u(t)}_{\infty} 
	+ |t|^{-\frac{d}{2}(1- \lambda) - \alpha(1- 2\lambda)} \( \norm{u(t)}_{2} + \norm{| J(t)|^{\beta} u(t)}_{2} \)
\end{align*}
for any $t \geq r_{0}$.
\end{Lem}
\begin{proof}
The proof is similar to that of \cite[Lemma 2.1]{HN98} and \cite[Lemma 2.11]{KM21} and hence it is omitted.
\end{proof}

We give the proof of Proposition \ref{prop4.5} dividing into two parts as follows:

\begin{Prop} \label{iv:pro1}
There exists $\varepsilon_1>0$ such that the following assertion holds: 
If $\varepsilon$, $u_0$, and $u$ satisfy the assumption in Proposition \ref{prop4.5}, 
then the estimate
\begin{align*}
	&{}|t|^{-A \varepsilon^{\frac{1}{d(1-\lambda)}}} \( \norm{u(t)}_{2} + \norm{|J(t)|^{\beta}u(t)}_{2} \) \leq \varepsilon 
\end{align*}
is valid for any $t \in I_{T}$.
\end{Prop}

\begin{proof}
Multiplying the integral equation associated with \eqref{eq1} with $u(r_{0}) = u_{0}$
by the operator $|J(t)|^{\beta}$, we see from Lemma \ref{lem:k1} that
\begin{align*}
	|J(t)|^{\beta} u(t) = U(t, r_{0}) |J(r_{0})|^{\beta} u_{0} + i \int_{r_{0}}^{t} U(t, s) |J(s)|^{\beta}F(u)(s))\, ds.
\end{align*}
We see from \eqref{eq4.15} and Lemma \ref{lem:34} that
\begin{align*}
	\norm{|J(t)|^{\beta} u(t)}_{2}
	\leq{}& \norm{U(t, r_{0}) |J(r_{0})|^{\beta} u_{0}}_{2} +\int_{r_{0}}^{t} \norm{|J(s)|^{\beta} F(u(s))}_{2} ds
	\\
	\leq{}& \norm{|J(r_{0})|^{\beta} u_{0}}_{2} + c_0 \varepsilon^{\frac{1}{d(1- \lambda)}} \int_{r_{0}}^{t} |s|^{-1} \norm{|J(s)|^{\beta} u(s)}_{2}\, ds.
\end{align*}
We here note $A = c_0$. 
Hence, the Gronwall inequality gives us 
\begin{align}
	|t|^{-A \varepsilon^{\frac{1}{d(1-\lambda)}}} \norm{|J(t)|^{\beta}u(t)}_{2} \leq \norm{|J(r_{0})|^{\beta} u_{0}}_{2}
	\label{iv:5}
\end{align}
for any $t \in I_{T}$.
Arguing as in the above, since 
\begin{align*}
	\norm{u(t)}_{2} 
	\leq{}& \norm{u_{0}}_{2} +  c_0  \varepsilon^{\frac{1}{d(1- \lambda)}} \int_{r_{0}}^{t} |s|^{-1} \norm{u(s)}_{2}\, ds,
\end{align*}
we get
\begin{align}
	|t|^{-A\varepsilon^{\frac{1}{d(1- \lambda)}}} \norm{u(t)}_{2} \leq \norm{u_{0}}_{2}  \label{iv:6}.
\end{align}
Collecting \eqref{iv:5} and \eqref{iv:6}, the desired estimate is valid. 
\end{proof}

\begin{Prop} \label{iv:pro2} 
There exists $\varepsilon_1 \in (0,1)$ such that the following assertion holds: 
If $\varepsilon$, $u_0$, and $u$ satisfy the assumption in Proposition \ref{prop4.5}, 
then the estimate
\begin{align*}
	&{}|t|^{\frac{d}{2}(1- \lambda)} \norm{u(t)}_{\infty} \leq C\varepsilon 
\end{align*}
is valid for any $t \in I_{T}$.
\end{Prop}

\begin{proof}
Take $\varepsilon_1 < \min\{1, \alpha (1- 2\lambda)((1+p_{c})A)^{-1}\}$ for some $\alpha \in (0, (\beta - d/2)/2)$. 
Combining Lemma \ref{iv:3} with Proposition \ref{iv:pro1}, we obtain
\begin{align}
	\norm{u(t)}_{\infty} 
	\leq{}& C|t|^{-\frac{d}{2}(1- \lambda)} \norm{\mathcal{F} \mathcal{M}_{+} U(0, t) u(t)}_{\infty} 
	+ C\varepsilon|t|^{-\frac{d}{2}(1- \lambda)- \alpha (1 -2\lambda)+A \varepsilon^{\frac{1}{d(1- \lambda)}}} 
	\label{iv:8}
\end{align}
for any $t \in I_{T}$.
Let us estimate the first term in the above last line.
Let $v$ satisfy $U(t,0) \mathcal{M}_{+}^{-1} v = u$, namely, $v = \mathcal{M}_{+} U(0, t)u$. 
Let us here introduce a phase modification
\begin{align*}
	 w_{p} \coloneqq \mathcal{F}^{-1}( \widetilde{\mathcal{P}}(t) \mathcal{F}  v), \qquad
	 \widetilde{\mathcal{P}}(t) \coloneqq \exp \left( i \eta \int_{r_{0}}^{t} |\mathcal{F}  v(\tau)|^{p_{c}} |\zeta_{2}(\tau)|^{- \frac{1}{1- \lambda}}\, d\tau  \right).
\end{align*}
Arguing as in the derivation of (3.14) in \cite{HN98}, 
we obtain
\begin{align*}
	i \partial_t (\mathcal{F} w_{p}) = \eta |\zeta_{2}(t)|^{-\frac{1}{1- \lambda}} \widetilde{\mathcal{P}}(t) \(I_{1}(t) + I_{2}(t) \),
\end{align*}
where $I_{j}(t)$ are defined by \eqref{remain:1}.
Integrating the above on $[r_{0}, t]$, noting that $\mathcal{F}\widetilde{w}(r_{0}) = \mathcal{F}v(r_{0}) = \mathcal{F} \mathcal{M}_{+} U(0, r_{0})u_{0}$, one has
\begin{align}
	\mathcal{F} w_{p}(t) = \mathcal{F} \mathcal{M}_{+} U(0, r_{0})u_{0} 
	- i \eta \int_{r_{0}}^{t} \widetilde{\mathcal{P}}(\tau) \( I_1(\tau)+I_2(\tau) \) |\zeta_{2}(\tau)|^{- \frac{1}{1- \lambda}}\, d \tau. \label{iv:12a}
\end{align}
Arguing as in \eqref{96:a2} and \eqref{iv:11}, one has
\begin{align}
	\begin{aligned}
	\norm{I_{1}(t)}_{\infty} + \norm{I_{2}(t)}_{\infty}
	\lesssim{}& |t|^{-\alpha (1 - 2 \lambda)} \( \norm{u(t)}_{2} + \norm{|J(t)|^{\beta} u(t)}_{2} \)^{p_{c}+1}
	\end{aligned}
	\label{iv:10}
\end{align}
for any $\alpha \in [0, ( \beta - d/2)/2)$. 
Moreover, by the assumption, we have
\begin{align*}
	 \norm{\mathcal{F} \mathcal{M}_{+} U(0, r_{0})u(r_{0})}_{\infty} \leq \norm{U(0, r_{0})u(r_{0}) }_{1} 
	 \leq C \(\norm{u_{0}}_{2} + \norm{|J(r_{0})|^{\beta} u_{0}}_{2} \) \leq C \varepsilon. 
\end{align*} 
Combining \eqref{iv:10} with Proposition \ref{iv:pro1}, we reach to
\begin{align}
	\begin{aligned}
	\norm{\mathcal{F} \mathcal{M}_{+} U(0, t)u(t)}_{\infty} 
	\leq{}&C\varepsilon+ C\varepsilon^{p_{c}+1} \int_{r_{0}}^{t} \tau^{-\alpha(1- 2\lambda)-1+ (p_{c}+1) A \varepsilon^{\frac{1}{d(1-\lambda)}}}\, d\tau \\
	\leq{}& C(\varepsilon_1, \alpha, \lambda) \varepsilon
	\end{aligned}
	\label{iv:9}
\end{align}
for any $t \in T_{T}$. 
Therefore collecting \eqref{iv:8} and \eqref{iv:9}, we conclude
\begin{align*}
	|t|^{\frac{d}{2}(1- \lambda)} \norm{u(t)}_{\infty} \leq C ( 1+ t^{-\alpha(1-2\lambda) + A\varepsilon^{\frac{1}{d(1-\lambda)}}})\varepsilon \leq C\varepsilon
\end{align*}
for any $t \in I_{T}$. This completes the proof. 
\end{proof}

\begin{proof}[\textit{\bf Proof of Proposition \ref{prop4.5}}]
The consequence immediately follows from Proposition \ref{iv:pro2} and Proposition \ref{iv:pro1}.
\end{proof}

\section{Derivation of the integral equation \eqref{inte:1}} \label{app:B}

In this appendix, we derive \eqref{inte:1}.
Let us consider the term $ i \partial_t \(U(0,t) u_{p}\)  - U(0,t) F(u_{p})$ in \eqref{intee}.
Note that $i \partial_{t} \widehat{w} = (c_{+} |s|)^{-1} F\( \widehat{w}\)$.
By using \eqref{mdfm}, we have
\begin{align*}
	U(0,t) F(u_{p}) 
	={}& \mathcal{M}_{2}^{-1} \mathcal{M}_{+}^{-1} \mathcal{F}^{-1} | \zeta_{2}(t)|^{- \frac{1}{1- \lambda}} F\( \widehat{w} \)  \\
	&{}+ \mathcal{M}_{+}^{-1} \( \mathcal{M}_{2}^{-1} -1 \) \mathcal{F}^{-1} (c_{+} |t|)^{-1} F\( \widehat{w} \) 
	+ \mathcal{M}_{+}^{-1} \mathcal{F}^{-1} i \partial_{t} \widehat{w}.
\end{align*}
A direct computation shows
\begin{align*}
	i \partial_t \(U(0,t) u_{p}\) = \mathcal{M}_{+}^{-1} \( \mathcal{M}_{2}^{-1} -1 \) \mathcal{F}^{-1} i \partial_{t} \widehat{w}
	+ \mathcal{M}_{+}^{-1} \mathcal{F}^{-1} i \partial_{t} \widehat{w}.
\end{align*}
Combining these above, we obtain
\begin{align*}
	&{}i \partial_t \(U(0,t) u_{p}(t)\)  - U(0,t) F(u_{p}(t)) \\
	={}& -\mathcal{M}_{+}^{-1} \mathcal{M}_{2}^{-1}(t) \mathcal{F}^{-1} \( |\zeta_{2}(t)|^{- \frac{1}{1- \lambda}} - (c_{+} |t|)^{-1} \) F\( \widehat{w}(t) \) \\
	&{}+ \mathcal{M}_{+}^{-1} \( 1 - \mathcal{M}_{2}^{-1}(t) \) \mathcal{F}^{-1} (c_{+} |t|)^{-1} F\( \widehat{w}(t) \) \\
	&{}- \mathcal{M}_{+}^{-1} \( 1 - \mathcal{M}_{2}^{-1}(t) \) \mathcal{F}^{-1} i \partial_{t} \widehat{w}(t).
\end{align*}
Since we may assume that $\norm{v(t)}_{2} \rightarrow 0$ as $t \rightarrow -\infty$, 
substituting the above into \eqref{intee} and integrating with respect to $t$, one sees from \eqref{mdfm} that
\begin{align*}
	v(t) 
	={}& -i \int_{- \infty}^{t} U(t,s) \{ F(u_{p}(s)+v(s)) - F(u_{p}(s)) \}\, ds \\
	&{}- i \int_{- \infty}^{t} U(t,0) \mathcal{M}_{+}^{-1} \mathcal{M}_{2}^{-1}(s) \mathcal{F}^{-1} \( |\zeta_{2}(s)|^{- \frac{1}{1- \lambda}} - (c_{+} |s|)^{-1} \) F\( \widehat{w} \)\, ds \\
	&{}+ i \int_{- \infty}^{t} U(t,s) \mathcal{M}_{1}(s) \mathcal{D}_{1}(s) \(\mathcal{F} \mathcal{M}_{2}(s) \mathcal{F}^{-1} -1 \) F\( \widehat{w}(s) \)\, \frac{ds}{c_{+} |s|} \\
	&{}+ \mathcal{M}_{1}(t) \mathcal{D}_{1}(t) \(\mathcal{F} \mathcal{M}_{2}(t) \mathcal{F}^{-1} -1 \) \widehat{w}(t)
\end{align*}
and hence this is \eqref{inte:1}.

\section*{Acknowledgments}
M.K. was supported by JSPS KAKENHI Grant Numbers 20K14328.
H.M. was supported by JSPS KAKENHI Grant Numbers 22K13941.
The authors express their gratitude to Professor Satoshi Masaki for his valuable suggestion regarding the Dollard type decomposition \eqref{mdfm}. 
The authors also would like to thank Professors Kota Uriya and Yuta Wakasugi for their helpful advice on the proof of Theorem \ref{thm:1}.

\bibliographystyle{amsplain}

\end{document}